\documentclass[11pt]{amsart}
\usepackage{mathrsfs}
\usepackage{amssymb}

\usepackage{amscd}
\usepackage{amsmath, amssymb}
\usepackage{amsfonts}
\usepackage[colorlinks,linkcolor=blue,citecolor=blue, pdfstartview=FitH]{hyperref}
\usepackage[all]{xy}

\numberwithin{equation}{section}

\theoremstyle{plain}
\newtheorem{thm}{Theorem}[section]
\newtheorem{theorem}[thm]{Theorem}
\newtheorem{lemma}[thm]{Lemma}
\newtheorem{corollary}[thm]{Corollary}
\newtheorem{proposition}[thm]{Proposition}

\theoremstyle{definition}

\newtheorem{question}[thm]{Question}

\newtheorem{remark}[thm]{Remark}

\newtheorem{definition}[thm]{Definition}

\newtheorem{example}[thm]{Example}

\newtheorem{defn-thm}[thm]{Definition-Theorem}

\newcommand{\Image}{{ \textrm{Im}~}}

\begin{document}

\title{On the deformed Bott-Chern cohomology}
\author{Wei Xia}
\address{Wei Xia, Mathematical Science Research Center, Chongqing University of Technology, Chongqing, P.R.China, 400054.} \email{xiaweiwei3@126.com}

\thanks{This work was supported by the National Natural Science Foundation of China No. 11901590.}
\date{\today}

\begin{abstract}
Given a compact complex manifold $X$ and a integrable Beltrami differential $\phi\in A^{0,1}(X, T_{X}^{1,0})$, we introduce a double complex structure on $A^{\bullet,\bullet}(X)$ naturally determined by $\phi$ and study its Bott-Chern cohomology. In particular, we establish a deformation theory for Bott-Chern cohomology and use it to compute the deformed Bott-Chern cohomology for the Iwasawa manifold and the holomorphically parallelizable Nakamura manifold. The $\partial\bar{\partial}_{\phi}$-lemma is studied and we show a compact complex manifold satisfying $\partial\bar{\partial}_{\phi}$-lemma is formal.

\vskip10pt
\noindent
{\bf MSC~Classification (2010):} 32G05, 32C35, 32G99
\end{abstract}
\maketitle

\section{Introduction}
The Bott-Chern cohomology are important invariants of complex manifolds~\cite{BC65}. It has been studied by many authors in recent years~\cite{Ang13,AT15a,AT15b,AT17,ADT16,AK17a}. For example, Schweitzer studied the Hodge theory for Bott-Chern cohomology and gave a hypercohomology interpretation to it~\cite{Sch07}. Angella-Tomassini proved Fr\"ohlicher type inequalities for Bott-Chern cohomology and gave a beautiful characterization of the $\partial\bar{\partial}$-lemma~\cite{AT13}. Recently, S. Yang and X. Yang proved a blow-up formula for the Bott-Chern cohomology and they showed that satisfying the $\partial\bar{\partial}$-Lemma is a bimeromorphic invariant for threefolds ~\cite{YY17}, see \cite{RYY19,ASTT17,Ste18a,Ste18b,Meng19} for related works.

Let $X$ be a complex manifold and $X_t$ a small deformation (of $X$) whose complex structure is represented by a Beltrami differential $\phi\in A^{0,1}(X, T_{X}^{1,0})$. In this paper, we will study the Bott-Chern cohomology of the double complex $(A^{\bullet,\bullet}(X), \partial, \bar{\partial}_{\phi})$:
\begin{equation}
H^{p,q}_{BC\phi}(X) := \frac{\ker d_{\phi} \cap A^{p,q}(X)}{\Image \partial\bar{\partial}_{\phi} \cap A^{p,q}(X)},
\end{equation}
which we called the \emph{deformed Bott-Chern cohomology}, where $d_{\phi}=\partial+\bar{\partial}_{\phi}$ and $\bar{\partial}_{\phi} = \bar{\partial} - \mathcal{L}_{\phi}^{1,0}$. In Section \ref{hypercohomology interpretation}, we will show that there are similar hypercohomology interpretations to the deformed Bott-Chern cohomology as to the usual Bott-Chern cohomology.

Let $\pi: (\mathcal{X}, X)\to (B,0)$ be a deformation of a compact complex manifold $X$ such that for each $t\in B$ the complex structure on $X_t$ is represented by Beltrami differential $\phi(t)$. Given a Bott-Chern class $[y]\in H_{BC}^{p,q}(X)$, as motivated by our previous work on deformation of Dolbeault cohomology classes~\cite{Xia19dDol}, we try to construct a family of $(p,q)$-forms $\sigma (t)$ (on an analytic subset $T$ of $B$) such that
\begin{itemize}
  \item[1.] $\sigma (t)$ is holomorphic in $t$;
  \item[2.] $\partial\sigma (t) =\bar{\partial}_{\phi(t)}\sigma (t)= 0,~\forall t\in T$;
  \item[3.] $[\sigma (0)] = [y]\in H_{BC}^{p,q}(X)$.
\end{itemize}
We will develop a deformation theory for Bott-Chern cohomology in this respect, see Section~\ref{Deformations of BC c}. Among other things, we show the following
\begin{theorem}[=Theorem \ref{th-analytic subset}]\label{th-analytic subset 0}
Let $\pi: (\mathcal{X}, X)\to (B,0)$ be a deformation of a compact complex manifold $X$ such that for each $t\in B$ the complex structure on $X_t$ is represented by Beltrami differential $\phi(t)$. Then the set $\{t\in B\mid \dim H_{BC\phi(t)}^{p,q}(X)\geq k\}$ is an analytic subset of $B$ for any nonnegative integer $k$.
\end{theorem}

In \cite[Thm.\,1\,and\,2]{AT15b}, Angella-Tomassini generalized their previous result~\cite{AT13} to arbitrary double complex~\cite{AT15a}. This result, when applied to our situation, will give rise to the following
\begin{theorem}[=Theorem \ref{thm-AT-Frolicher-inequality-deformed}]\label{thm-AT-Frolicher-inequality-deformed 0}
Let $X$ be a compact complex manifold and $X_t$ a small deformation (of $X$) whose complex structure is represented by a Beltrami differential $\phi\in A^{0,1}(X, T_{X}^{1,0})$. Then for every $(p,q)\in \mathbb{N}\times \mathbb{N}$, we have
\begin{equation}\label{eq-AT-Frolicher-inequality-deformed 0}
\dim H_{BC\phi}^{p,q}(X)+ \dim H_{A\phi}^{p,q}(X) \geq \dim H_{\bar{\partial}_t}^{p,q}(X_t)+ \dim H_{\partial}^{p,q}(X).
\end{equation}
In particular, for every $k\in \mathbb{N}$, we have
\begin{equation}\label{eq-ddbar-phi-lemma-kehua 0}
\sum_{p+q=k}\dim H_{BC\phi}^{p,q}(X)+ \sum_{p+q=k}\dim H_{A\phi}^{p,q}(X) \geq 2\dim H_{dR}^{k}(X),
\end{equation}
and equality holds if and only if $X$ satisfies the $\partial\bar{\partial}_{\phi}$-lemma.
\end{theorem}

Note that when $X_t$ is a trivial deformation, i.e. $\phi=0$, Theorem \ref{eq-AT-Frolicher-inequality-deformed 0} is reduced to the result in ~\cite{AT13}. Combine Theorem \ref{eq-AT-Frolicher-inequality-deformed 0} with Theorem \ref{th-analytic subset 0}, we get
\begin{corollary}\label{coro-deformation-closedness-ddbar-lemma 0}
Let $\pi: (\mathcal{X}, X)\to (B,0)$ be a small deformation of the compact complex manifold $X$ such that for each $t\in B$ the complex structure on $X_t$ is represented by Beltrami differential $\phi(t)$. Then
the set
\[
T:=\{t\in B\mid X~\text{satisfies the}~\partial\bar{\partial}_{\phi(t)}\text{-lemma} \}
\]
is an analytic open subset (i.e. complement of analytic subset) of $B$. In particular, if $B\subset \mathbb{C}$ is a small open disc with $0\in B$ and $T$ is not empty, then $T=B$ or $T=B\setminus\{0\}$.
\end{corollary}
It is known that satisfying the $\partial\bar{\partial}$-lemma is a deformation open property and not a deformation closed property in the sense of Popovici~\cite{Pop14}, see ~\cite{Wu06,AT13,AK17b} and the references therein. But it is still not clear whether satisfying the $\partial\bar{\partial}$-lemma is an analytically open property, i.e. does the corresponding statement in Corollary \ref{coro-deformation-closedness-ddbar-lemma 0} holds for the $\partial\bar{\partial}$-lemma? On the other hand, we see from Corollary \ref{coro-deformation-closedness-ddbar-lemma 0} that if $X$ satisfies the $\partial\bar{\partial}$-lemma then $X$ also satisfies the $\partial\bar{\partial}_{\phi(t)}$-lemma for small $t$. But conversely, if $X$ satisfies the $\partial\bar{\partial}_{\phi(t)}$-lemma for all small $t\neq 0$ it is possible that $X$ does not satisfy the $\partial\bar{\partial}$-lemma\footnote{Though it is still not known whether such examples exist, we think they should be large in number.}. Hence the following Theorem generalize the corresponding well-known result of Deligne-Griffiths-Morgan-Sullivan~\cite{DGMS75}:
\begin{theorem}\label{thm-formal 0}
Let $X$ be a compact complex manifold and $X_t$ a small deformation (of $X$) whose complex structure is represented by a Beltrami differential $\phi\in A^{0,1}(X, T_{X}^{1,0})$. If $X$ satisfies the $\partial\bar{\partial}_\phi$-lemma, then $X$ is formal.
\end{theorem}

The dimensions of the deformed Bott-Chern cohomology is computed for the Iwasawa manifold and the holomorphically parallelizable Nakamura manifold, see Section \ref{The deformed Bott-Chern cohomology of the Iwasawa manifold and the holomorphically parallelizable Nakamura manifold}. Comparing this with the computations of Angella-Kasuya~\cite{AK17b}, we see that there exists compact complex manifold $X$ and its small deformation $X_t$ such that $X_t$ satisfy the $\partial\bar{\partial}$-lemma but $X$ does not satisfy the $\partial\bar{\partial}_{\phi(t)}$-lemma.

There are many questions regarding the $\partial\bar{\partial}_{\phi}$-lemma may be asked:
\begin{question}Let $\pi: (\mathcal{X}, X)\to (B,0)$ be a small deformation of the compact complex manifold $X$ such that for each $t\in B$ the complex structure on $X_t$ is represented by Beltrami differential $\phi(t)$.
\begin{itemize}
  \item[1.] Is it true that
\begin{equation}\label{eq-dim-comparison}
\dim H_{BC\phi(t)}^{p,q}(X)\geq\dim H_{BC}^{p,q}(X_t)
 \end{equation}
for any $t\in B$ and $(p,q)\in \mathbb{N}\times \mathbb{N}$? If this holds, then $X$ satisfies the $\partial\bar{\partial}_{\phi(t)}$-lemma will imply $X_t$ satisfy the $\partial\bar{\partial}$-lemma. Note that \eqref{eq-dim-comparison} is true for the examples considered in Section \ref{The deformed Bott-Chern cohomology of the Iwasawa manifold and the holomorphically parallelizable Nakamura manifold};
  \item[2.] If $B\subset \mathbb{C}$ is a small open disc with $0\in B$, can we find an example such that $T=B\setminus\{0\}$ (in the notation of Corollary \ref{coro-deformation-closedness-ddbar-lemma 0})? According to Corollary \ref{coro-deformation-closedness-ddbar-lemma 0}, there should be many such examples. In this case, the Fr\"ohlicher spectral sequence on the central fiber $X$ must degenerates at $E_1$, see Remark~\ref{rk-to-AK-Frolicher-ineq};
  \item[3.] If $X_t$ is K\"ahler, is it true that $X$ must satisfy the $\partial\bar{\partial}_{\phi(t)}$-lemma?
\end{itemize}
\end{question}

\section{The deformed double complex $(A^{\bullet,\bullet}(X), \partial, \bar{\partial}_{\phi})$ and its Bott-Chern cohomology}
\label{The deformed double complex and its Bott-Chern cohomology}
Let $X$ be a complex manifold and $X_t$ a small deformation (of $X$) whose complex structure is represented by a Beltrami differential $\phi\in A^{0,1}(X, T_{X}^{1,0})$. Recall the following useful facts~\cite{LRY15,Xia19deri}:
\[
e^{-i_{\phi}}de^{i_{\phi}}=d-\mathcal{L}_{\phi}^{1,0}-\mathcal{L}_{\phi}^{0,1} -i_{\frac{1}{2}[\phi,\phi]}~\text{and}~\mathcal{L}_{\phi}^{0,1}=-i_{\bar{\partial}\phi}~.
\]
Since $\phi$ satisfy the Maurer-Cartan equation $\bar{\partial}\phi-\frac{1}{2}[\phi,\phi]=0$, we have
\begin{equation}\label{}
d_{\phi}~: = e^{-i_{\phi}}de^{i_{\phi}} = \partial + \bar{\partial}_{\phi},~~ \text{with}~~ \bar{\partial}_{\phi} = \bar{\partial} - \mathcal{L}_{\phi}^{1,0},
\end{equation}
and
\begin{equation}\label{}
d_{\bar{\phi}}~: = e^{-i_{\bar{\phi}}}de^{i_{\bar{\phi}}} = \partial_{\phi} + \bar{\partial},~~ \text{with}~~ \partial_{\phi} = \partial - \mathcal{L}_{\bar{\phi}}^{0,1}.
\end{equation}
Since $[\partial,\bar{\partial}_{\phi}]=[\partial_{\phi},\bar{\partial}]=0$, the \emph{deformed Bott-Chern cohomology} can be defined as follows:
\begin{equation}\label{}
H^{p,q}_{BC\phi}(X) := \frac{\ker d_{\phi} \cap A^{p,q}(X)}{\Image \partial\bar{\partial}_{\phi} \cap A^{p,q}(X)},~~~~~~
H^{p,q}_{BC\bar{\phi}}(X) := \frac{\ker d_{\bar{\phi}} \cap A^{p,q}(X)}{\Image \partial_{\phi}\bar{\partial} \cap A^{p,q}(X)},~~\forall p,q\geq 0~,
\end{equation}
and $h^{p,q}_{BC\phi}:=\dim H^{p,q}_{BC\phi}(X), h^{p,q}_{BC\bar{\phi}}:=\dim H^{p,q}_{BC\bar{\phi}}(X)$. The conjugation gives a natural isomorphism between $H^{p,q}_{BC\phi}(X)$ and $H^{q,p}_{BC\bar{\phi}}(X)$, we thus have $h^{p,q}_{BC\phi}=h^{q,p}_{BC\bar{\phi}}$.

\section{Hypercohomology interpretations to the deformed Bott-Chern cohomology}
\label{hypercohomology interpretation}
It is clear that the Poincar\'e lemma holds for $d_{\phi}$ and $\bar{\partial}_{\phi}$ (for the latter, see~\cite[Thm.\,3.4]{Xia19dDol}). The sheaf of germs of $\bar{\partial}_{\phi}$-closed $p$-forms will be denoted by $\Omega^{p}_{\phi}$. The following Lemma is essentially proved in~\cite{Sch07}:
\begin{lemma}
\label{lem-local solvability}
Let $U\subset\mathbb{C}^n$ be an open ball.
\begin{itemize}
  \item[1.] Let $\theta\in A^k(U)$ with $k\geq 1$ such that $\theta^{p,q}=0$ except $p_1\leq p\leq p_2(p_1< p_2)$. If $\theta$ is $d_{\phi}$-closed, then $\theta= d_{\phi} \alpha$ for some $\alpha\in A^{k-1}(U)$ with $\alpha^{p,q}=0$ except $p_1\leq p\leq p_2-1$.

  \item[2.] Assume $\theta \in A^{p,q}(U)$ is $d_{\phi}$-closed.
  \item[$i)$] If $p\geq 1$ and $q\geq 1$, then $\theta\in \partial\bar{\partial}_{\phi} A^{p-1,q-1}(U)$.
  \item[$ii)$] If $p\geq 1$ and $q=0$, then $\theta\in \partial \Omega^{p-1}_{\phi}(U)$.
  \item[$\bar{ii})$] If $p=0$ and $q\geq 1$, then $\theta\in \bar{\partial}_{\phi} \bar{\Omega}^{q-1}(U)$.
  \item[$iii)$] If $p=q=0$, then $\theta$ is a constant.

  \item[3.] Assume $\theta \in A^{p,q}(U)$ is $\partial\bar{\partial}_{\phi}$-closed.
  \item[$i)$] If $p\geq 1$ and $q\geq 1$, then $\theta\in \bar{\partial}_{\phi} A^{p,q-1}(U)+\partial A^{p-1,q}(U)$.
  \item[$ii)$] If $p\geq 1$ and $q=0$, then $\theta\in \Omega_{\phi}^{p}(U) + \partial A^{p-1,0}(U)$.
  \item[$\bar{ii})$] If $p=0$ and $q\geq 1$, then $\theta\in \bar{\partial}_{\phi}A^{0,q-1}(U) + \bar{\Omega}^{q}(U)$.
  \item[$iii)$] If $p=q=0$, then $\theta\in \mathcal{O}_{\phi}(U) + \bar{\mathcal{O}}(U)$.

  \item[4.] Let $\theta\in A^k(U)$ with $k\geq 1$ and $p_1,~q_1,~p_2,~q_2$ be two positive integers with $p_1+q_1=p_2+q_2= k$. If $(d_{\phi}\theta)^{p,q}=0$ for $p+q=k+1,~p_1+1\leq p\leq p_2$ and $q_1\geq q\geq q_2+1$, then there exists $\gamma^{p_1,q_1}, \alpha^{p_1,q_1-1}, \alpha^{p_1+1,q_1-2},\cdots, \alpha^{p_2-1,q_2}, \gamma^{p_2,q_2}$, s.t. $\gamma^{p_1,q_1}$ is $\partial$-closed, $\gamma^{p_2,q_2}$ is $\bar{\partial}_{\phi}$-closed and
\begin{align*}
\theta^{p_1,q_1}&=\gamma^{p_1,q_1}+\bar{\partial}_{\phi}\alpha^{p_1,q_1-1},\\
\theta^{p_1+1,q_1-1}&=\partial\alpha^{p_1,q_1-1}+\bar{\partial}_{\phi}\alpha^{p_1+1,q_1-2},\\
 &\cdots,\\
\theta^{p_2-1,q_2+1}&=\partial\alpha^{p_2-2,q_2+1}+\bar{\partial}_{\phi}\alpha^{p_2-1,q_2},\\
\theta^{p_2,q_2}&=\partial\alpha^{p_2-1,q_2}+\gamma^{p_2,q_2},
\end{align*}
in particular, we have
\[
\theta^{p_1,q_1} + \theta^{p_1+1,q_1-1}+ \cdots + \theta^{p_2,q_2} = \gamma^{p_1,q_1} + d_{\phi}\alpha+\gamma^{p_2,q_2},
\]
where $\alpha=\alpha^{p_1,q_1-1}+\alpha^{p_1+1,q_1-2}+\cdots +\alpha^{p_2-1,q_2}$.
\end{itemize}
\end{lemma}

\begin{proof}1.First, by the $d_{\phi}$-Poincar\'e lemma, we can write $\theta= d_{\phi}\beta$ for some $\beta\in A^{k-1}(U)$. If $p_1=0$ and $p_2=k$ there is nothing to prove, so we assume\footnote{We may further assume that $k\geq 2 $ because the case $k= 1$ is trivial.} $p_1>0$ or $p_2<k$. We first consider the case $p_1>0$. We deduce from $\theta= d_{\phi}\beta$ that $\bar{\partial}_{\phi}\beta^{0,k-1}= \theta^{0,k}=0$, and by applying the $\bar{\partial}_{\phi}$-Poincar\'e lemma, one can write $\beta^{0,k-1}= \bar{\partial}_{\phi}\gamma^{0,k-2}$. Set $\tilde{\beta}:=\beta-d_{\phi}\gamma^{0,k-2}$, we have $d_{\phi}\tilde{\beta}=\theta$ but $\tilde{\beta}^{0,k-1}=0$. We can therefore assume that $\beta$ does not have components of type $(0,k-1)$. Now if $p_1>1$, then since $\beta^{0,k-1}=0$ we have $0= \theta^{1,k-1}=\bar{\partial}_{\phi}\beta^{1,k-2}+\partial\beta^{0,k-1}=\bar{\partial}_{\phi}\beta^{1,k-2}$. By the $\bar{\partial}_{\phi}$-Poincar\'e lemma, one can write $\beta^{1,k-2}= \bar{\partial}_{\phi}\gamma^{0,k-3}$. Set $\tilde{\beta}:=\beta-d_{\phi}\gamma^{0,k-3}$, we have $d_{\phi}\tilde{\beta}=\theta$ but $\tilde{\beta}^{1,k-2}=0$. We can therefore assume that $\beta$ does not have components of type $(1,k-2)$. By repeating this reasoning, we can assume that $\beta$ does not have components of type $(p,q)$ for $p<p_1$. The case $p_2<k$ can be proved in the same way by applying the $\partial$-Poincar\'e lemma.

2.$iii)$ is obvious. We first assume $p\geq 1$. We apply 1. to the form $\theta$ for $p_1=p-1, p_2=p$: there exists $\alpha\in A^{p-1,q}(U)$ s.t. $\theta= d_{\phi}\alpha$ and so $\theta= \partial\alpha$ with $\bar{\partial}_{\phi}\alpha=0$. This is $ii)$. If furthermore $q\geq 1$, by the $\bar{\partial}_{\phi}$-Poincar\'e lemma, we can write $\alpha=\bar{\partial}_{\phi}\beta$ and so $\theta= \partial\bar{\partial}_{\phi}\beta$. This is $i)$. For $\bar{ii})$, we apply 1. to $\theta$ for $p_1=0, p_2=1$: there exists $\alpha\in A^{0,q-1}(U)$ s.t. $\theta= d_{\phi}\alpha$ and so $\theta= \bar{\partial}_{\phi}\alpha$ with $\partial\alpha=0$.

3.Set $\theta^{p+1,q}:= \partial\alpha^{p,q}$ then $\theta^{p+1,q}$ is $d_{\phi}$-closed. By 2.$i)$ and $ii)$, there exists $\alpha\in A^{p,q}(U)$ s.t. $\theta^{p+1,q}= \partial\alpha$ with $\bar{\partial}_{\phi}\alpha=0$. Note that $\partial(\theta-\alpha)=0$ and $\theta=(\theta-\alpha)+ \alpha$. Then 3. follows from the $\bar{\partial}_{\phi}$-Poincar\'e lemma and the $\partial$-Poincar\'e lemma.

4.First from the assumption we see that $(d_{\phi}\theta)^{p_1+1,q_1}=\partial\theta^{p_1,q_1}+ \bar{\partial}_{\phi} \theta^{p_1+1,q_1-1}=0$. In particular, $\theta^{p_1,q_1}$ is $\partial\bar{\partial}_{\phi}$-closed. By 3.$i)$ and $\bar{ii})$ there exists $\gamma^{p_1,q_1}$ s.t. $\gamma^{p_1,q_1}$ is $\partial$-closed\footnote{$\gamma^{p_1,q_1}$ is $\partial$-exact if $p_1\geq 1$.} and $\theta^{p_1,q_1}=\gamma^{p_1,q_1}+\bar{\partial}_{\phi}\alpha^{p_1,q_1-1}$. Note that $\bar{\partial}_{\phi}( \theta^{p_1+1,q_1-1}-\partial\alpha^{p_1,q_1-1})=\partial (-\theta^{p_1,q_1}+\bar{\partial}_{\phi}\alpha^{p_1,q_1-1})=\partial\gamma^{p_1,q_1}=0$, so $\theta^{p_1+1,q_1-1}=\partial\alpha^{p_1,q_1-1}+\gamma^{p_1+1,q_1-1}$ such that $\gamma^{p_1+1,q_1-1}$ is $\bar{\partial}_{\phi}$-closed. Hence we can write $\theta^{p_1+1,q_1-1}=\partial\alpha^{p_1,q_1-1}+\bar{\partial}_{\phi}\alpha^{p_1+1,q_1-2}$.

Again from the assumption we see that $(d_{\phi}\theta)^{p_1+2,q_1-1}=\partial\theta^{p_1+1,q_1-1}+ \bar{\partial}_{\phi} \theta^{p_1+2,q_1-2}=0$, and note that $\bar{\partial}_{\phi}( \theta^{p_1+2,q_1-2}-\partial\alpha^{p_1+1,q_1-2})=\partial (-\theta^{p_1+1,q_1-1}+\bar{\partial}_{\phi}\alpha^{p_1+1,q_1-2})=-\partial^2\alpha^{p_1,q_1-1}=0$, we have $\theta^{p_1+2,q_1-2}=\partial\alpha^{p_1+1,q_1-2}+\bar{\partial}_{\phi}\alpha^{p_1+2,q_1-3}$.

Continuing in this way, we get the desired results. In the last two steps, from $(d_{\phi}\theta)^{p_2-1,q_2+2}=0$ we get $\theta^{p_2-1,q_2+1}=\partial\alpha^{p_2-2,q_2+1}+\bar{\partial}_{\phi}\alpha^{p_2-1,q_2}$ and from $(d_{\phi}\theta)^{p_2,q_2+1}=0$ we get $\theta^{p_2,q_2}=\partial\alpha^{p_2-1,q_2}+\gamma^{p_2,q_2}$.
\end{proof}

Let $X$ be a complex manifold. For fixed $p\geq 1$ and $q\geq 1$, we define a sheaf complex $\mathscr{L}_{\phi}^\bullet$ (which depend on $(p,q)$) as follows:
\begin{equation}
\left\{
\begin{array}{ll}
\mathscr{L}_{\phi}^k= \bigoplus_{r+s=k,r<p,s<q}A^{r,s},~~~\text{for}~~~k\leq p+q-2, & \\
\mathscr{L}_{\phi}^{k-1}= \bigoplus_{r+s=k,r\geq p,s\geq q}A^{r,s},~~~\text{for}~~~k\geq p+q.  & \\
\end{array} \right.
\end{equation}
The differential is given by
\begin{equation*}
\xymatrix@C=0.5cm{
  0 \ar[r] & \mathscr{L}_{\phi}^0 \ar[rr]^{\Pi_{\mathscr{L}_{\phi}^1}d_{\phi}} && \mathscr{L}_{\phi}^1 \ar[rr]^{\Pi_{\mathscr{L}_{\phi}^2}d_{\phi}} && \mathscr{L}_{\phi}^2 \ar[r] & \cdots \\  \ar[r] & \mathscr{L}_{\phi}^{p+q-3} \ar[rr]^{\Pi_{\mathscr{L}_{\phi}^{p+q-2}}d_{\phi}} && \mathscr{L}_{\phi}^{p+q-2} \ar[rr]^{\partial\bar{\partial}_{\phi}} && \mathscr{L}_{\phi}^{p+q-1} \ar[rr]^{d_{\phi}} && \mathscr{L}_{\phi}^{p+q} \ar[rr]^{d_{\phi}} && \cdots ,}
\end{equation*}
where $\Pi_{\mathscr{L}^k}:\bigoplus_{r+s=k}A^{r,s}\longrightarrow \mathscr{L}_{\phi}^k$ is the projection.
In particular, we find that
\begin{equation*}
\xymatrix@C=0.5cm{
\mathscr{L}_{\phi}^{p+q-2}=A^{p-1,q-1} \ar[rr]^{\quad\partial\bar{\partial}_{\phi}} && \mathscr{L}_{\phi}^{p+q-1}=A^{p,q} \ar[rr]^{d_{\phi}\qquad} && \mathscr{L}_{\phi}^{p+q}=A^{p,q+1}\oplus A^{p+1,q},}
\end{equation*}
and so $\mathbb{H}^{p+q-1}(X,\mathscr{L}_{\phi}^\bullet)\cong H^{p+q-1}(\mathscr{L}_{\phi}^\bullet(X))=H_{BC\phi}^{p,q}(X)$. The sheaf complex $\mathscr{L}_{\phi}^\bullet$ has the following subcomplexes
\begin{equation*}
\xymatrix@C=0.5cm{
(\mathscr{S'}_{\phi}^\bullet,\partial):~\mathcal{O}_{\phi}\ar[r]^{\qquad\partial} & \Omega^{1}_{\phi} \ar[r]^{\partial} & \Omega^{2}_{\phi} \ar[r]^{\partial} & \cdots  \ar[r]^{\partial} & \Omega^{p-1}_{\phi} \ar[r] & 0,\\
(\mathscr{S''_{\phi}}^\bullet,\bar{\partial}_{\phi}):~\bar{\mathcal{O}}\ar[r]^{\qquad\bar{\partial}_{\phi}} & \bar{\Omega}^{1} \ar[r]^{\bar{\partial}_{\phi}} & \bar{\Omega}^{2} \ar[r]^{\bar{\partial}_{\phi}} & \cdots \ar[r]^{\bar{\partial}_{\phi}} & \bar{\Omega}^{q-1} \ar[r] & 0,}
\end{equation*}
and
\[
\mathscr{S}_{\phi}^\bullet:=(\mathscr{S'}_{\phi}^\bullet,\partial)+(\mathscr{S''}_{\phi}^\bullet,\bar{\partial}_{\phi})\footnote{The sum is direct except $k=0$ and $\mathcal{O}_{\phi}+\bar{\mathcal{O}}\longrightarrow \Omega^{1}_{\phi}\oplus\bar{\Omega}^{1}: f+g\to (\partial f, \bar{\partial}_{\phi} g)$.}.
\]
Note that by Lemma \ref{lem-local solvability}, the complex $(\mathscr{S'}_{\phi}^\bullet,\partial)$ is exact for $0<k< p-1$ where $\mathscr{S'}_{\phi}^k=\Omega^{k}_{\phi}$.
\begin{proposition}\label{prop-quasiiso-S-L}
The inclusion $\mathscr{S}_{\phi}^\bullet\hookrightarrow \mathscr{L}_{\phi}^\bullet$ induces an isomorphism $\mathscr{H}^k(\mathscr{S}_{\phi}^\bullet)\cong\mathscr{H}^k(\mathscr{L}_{\phi}^\bullet),~\forall k\geq 0$, and we have\footnote{See also~\cite[pp.\,31]{Koo11}.}
\begin{equation*}
\mathscr{H}^k(\mathscr{S}_{\phi}^\bullet)\cong\mathscr{H}^k(\mathscr{L}_{\phi}^\bullet)=\left\{
\begin{array}{ll}
\mathbb{C},~~~\text{for}~~~k=0,~p>1,~q>1, & \\
\mathcal{O}_{\phi},~~~\text{for}~~~k=0,~p=1,~q>1, & \\
\bar{\mathcal{O}},~~~\text{for}~~~k=0,~p>1,~q=1, & \\
\mathcal{O}_{\phi}\oplus\bar{\mathcal{O}},~~~\text{for}~~~k=0,~p=1,~q=1, & \\
\Omega^{p-1}_{\phi}/\partial\Omega^{p-2}_{\phi},~~~\text{for}~~~0<k=p-1~\text{and}~p\neq q,  & \\
\bar{\Omega}^{q-1}/\bar{\partial}_{\phi}\bar{\Omega}^{q-2},~~~\text{for}~~~0<k=q-1~\text{and}~p\neq q,  & \\
\Omega^{p-1}_{\phi}/\partial\Omega^{p-2}_{\phi}\oplus \bar{\Omega}^{p-1}/\bar{\partial}_{\phi}\bar{\Omega}^{p-2},~~~\text{for}~~~0<k=p-1=q-1,  & \\
0,~~~~~~\text{otherwise}.
\end{array} \right.
\end{equation*}
\end{proposition}

\begin{proof}First, we show that $\mathscr{H}^k(\mathscr{L}_{\phi}^\bullet)=0$ for $k\geq \max\{p,q\}$. In fact, for $k\geq p+q$, this follows from Lemma \ref{lem-local solvability} 1.; for $k= p+q-1$, this follows from Lemma \ref{lem-local solvability} 2.; for $k= p+q-2$, this follows from Lemma \ref{lem-local solvability} 3.; for $k< p+q-2$, this follows\footnote{We apply Lemma \ref{lem-local solvability} 4. for $p_1=k-q+1,~q_1=q-1,~p_2=p-1,~q_2=k-p+1$. Note that we have $\theta^{k-q+1,q-1}=\partial\gamma^{k-q,q-1}+\bar{\partial}_{\phi}\alpha^{k-q+1,q-2}$ and $d_{\mathscr{L}^{k-1}}\gamma^{k-q,q-1}=\partial\gamma^{k-q,q-1}$, where $d_{\mathscr{L}^{k-1}}=\Pi_{\mathscr{L}^{k-1}}d_{\phi}$. Similarly, $\theta^{p-1,k-p+1}=\partial\alpha^{p-2,k-p+1}+\bar{\partial}_{\phi}\gamma^{p-1,k-p}$ and $d_{\mathscr{L}^{k-1}}\gamma^{p-1,k-p}=\bar{\partial}_{\phi}\gamma^{p-1,k-p}$.} from Lemma \ref{lem-local solvability} 4.~.

Now we discuss the cases when $k<p$ or $k<q$.

For $k=p-1\geq q$, if $\theta=\theta^{p-q,q-1}+\cdots \theta^{p-1,0}\in\mathscr{L}_{\phi}^{p-1}(U)$ is $d_{\mathscr{L}_{\phi}^{p-1}}$-closed where $U\subset X$ is an open ball. By Lemma \ref{lem-local solvability} 4., we can write
\[
\theta^{p-q,q-1}=\gamma^{p-q,q-1}+\bar{\partial}_{\phi}\alpha^{p-q,q-2},~\cdots~,\theta^{p-1,0}=\partial\alpha^{p-2,0}+\gamma^{p-1,0},
\]
where $\gamma^{p-q,q-1}$ is $\partial$-closed and $\gamma^{p-1,0}$ is $\bar{\partial}_{\phi}$-closed. Since $p-q\geq 1$, we have $\gamma^{p-q,q-1}=\partial \gamma^{p-q-1,q-1}=d_{\mathscr{L}^{p-2}}\gamma^{p-q-1,q-1}$ and so
\[
\theta=d_{\mathscr{L}_{\phi}^{p-2}}(\gamma^{p-q-1,q-1}+\alpha)+\gamma^{p-1,0},~\text{with}~\alpha=\alpha^{p-q,q-2}+\cdots+\alpha^{p-2,0}.
\]
On the other hand, if $\theta$ is $d_{\mathscr{L}_{\phi}^{\bullet}}$-exact, then there exists $u=u^{p-q-1,q-1}+\cdots+ u^{p-2,0}\in\mathscr{L}_{\phi}^{p-2}(U)$ s.t.
\[
d_{\mathscr{L}_{\phi}^{p-2}}u=(d_{\phi}u)^{p-q,q-1}+\cdots+ (d_{\phi}u)^{p-1,0}=\theta=\theta^{p-q,q-1}+\cdots+ \theta^{p-1,0}.
\]
Therefore $\partial u^{p-2,0}=\theta^{p-1,0}=\partial\alpha^{p-2,0}+\gamma^{p-1,0}\Rightarrow~\gamma^{p-1,0}=\partial (u^{p-2,0}-\alpha^{p-2,0})$ and $u^{p-2,0}-\alpha^{p-2,0}$ is $\partial\bar{\partial}_{\phi}$-closed. By Lemma \ref{lem-local solvability} 3.$ii)$, we see that $\gamma^{p-1,0}\in \partial\Omega_{\phi}^{p-2}(U)$. We thus have
\[
\mathscr{H}^{p-1}(\mathscr{L}_{\phi}^\bullet)=\frac{\Image d_{\mathscr{L}_{\phi}^{p-2}} +\Omega^{p-1}_{\phi} }{\Image d_{\mathscr{L}_{\phi}^{p-2}} +\partial\Omega^{p-2}_{\phi} }=\frac{\Omega^{p-1}_{\phi}}{\partial\Omega^{p-2}_{\phi}}=\mathscr{H}^{p-1}(\mathscr{S}_{\phi}^\bullet).
\]

For $k=p-1<q-1$, if $\theta=\theta^{0,p-1}+\cdots \theta^{p-1,0}\in\mathscr{L}_{\phi}^{p-1}(U)$ is $d_{\mathscr{L}_{\phi}^{p-1}}$-closed, by Lemma \ref{lem-local solvability} 4., we can write
\[
\theta^{0,p-1}=\gamma^{0,p-1}+\bar{\partial}_{\phi}\alpha^{0,p-2},~\cdots,~\theta^{p-1,0}=\partial\alpha^{p-2,0}+\gamma^{p-1,0},
\]
where $\gamma^{0,p-1}$ is $\partial$-closed and $\gamma^{p-1,0}$ is $\bar{\partial}_{\phi}$-closed. Note that since $k=p-1<q-1$, we have $d_{\mathscr{L}_{\phi}^{p-1}}\theta=0\Rightarrow\bar{\partial}_{\phi}\theta^{0,p-1}=(d_{\phi}\theta)^{0,p}=0\Rightarrow\gamma^{0,p-1}\in \bar{\partial}_{\phi}\bar{\Omega}^{p-2}$ by Lemma \ref{lem-local solvability} 2.$\bar{ii})$. Hence $\gamma^{0,p-1}\in \Image d_{\mathscr{L}_{\phi}^{p-2}}$. On the other hand, if $\theta$ is $d_{\mathscr{L}_{\phi}^{\bullet}}$-exact, then one can show as above that $\gamma^{p-1,0}\in \partial\Omega_{\phi}^{p-2}(U)$. We thus have
\[
\mathscr{H}^{p-1}(\mathscr{L}_{\phi}^\bullet)=\frac{\Image d_{\mathscr{L}_{\phi}^{p-2}} +\Omega^{p-1}_{\phi} }{\Image d_{\mathscr{L}_{\phi}^{p-2}} +\partial\Omega^{p-2}_{\phi} }=\frac{\Omega^{p-1}_{\phi}}{\partial\Omega^{p-2}_{\phi}}=\mathscr{H}^{p-1}(\mathscr{S}_{\phi}^\bullet).
\]

 For $k=p-1= q-1$, we have
\[
\mathscr{H}^{p-1}(\mathscr{L}_{\phi}^\bullet)=\frac{\Image d_{\mathscr{L}_{\phi}^{p-2}} +\bar{\Omega}^{p-1}+\Omega^{p-1}_{\phi} }{\Image d_{\mathscr{L}_{\phi}^{p-2}}+\bar{\partial}_{\phi}\bar{\Omega}^{p-2} +\partial\Omega^{p-2}_{\phi}}=\frac{\Omega^{p-1}_{\phi}}{\partial\Omega^{p-2}_{\phi}}\oplus\frac{\bar{\Omega}^{p-1} }{\bar{\partial}_{\phi}\bar{\Omega}^{p-2}}   =\mathscr{H}^{p-1}(\mathscr{S}_{\phi}^\bullet).
\]
\end{proof}

Consider the complex $\mathscr{B}_{\phi}^\bullet$ which is a modification of $\mathscr{S}_{\phi}^\bullet$ given by\footnote{This is the case when $p\geq q$, the case $p< q$ is similar. To make our notations clear and simple, we will only write explicitly one of the cases in what follows.}
\begin{equation*}
\xymatrix@C=0.5cm{
\mathscr{B}_{\phi}^\bullet:~\mathbb{C}\ar[r] &\mathcal{O}_{\phi}\oplus\bar{\mathcal{O}} \ar[r]^{\partial\oplus\bar{\partial}_{\phi}} & \Omega^{1}_{\phi}\oplus\bar{\Omega}^{1} \ar[r]^{\partial\oplus\bar{\partial}_{\phi}} & \Omega^{2}_{\phi}\oplus\bar{\Omega}^{2} \ar[r] & \cdots \\
 \ar[r] & \Omega^{q-1}_{\phi}\oplus\bar{\Omega}^{q-1}\ar[r]^{\qquad \partial\oplus 0} & \Omega^{q}_{\phi}\ar[r]^{\partial} &\cdots  \ar[r]^{\partial} & \Omega^{p-1}_{\phi} \ar[r] & 0,}
\end{equation*}
where the first morphism is defined by
\[
\mathbb{C}\longrightarrow\mathcal{O}_{\phi}\oplus\bar{\mathcal{O}}:a\mapsto (a,-a)~.
\]

\begin{proposition}\label{prop-quasiiso-B-S-phi}
The natural map from $\mathscr{B}_{\phi}^\bullet$ to $\mathscr{S}_{\phi}^\bullet[1]$, where
\[
\mathscr{B}_{\phi}^1=\mathcal{O}_{\phi}\oplus\bar{\mathcal{O}}\longrightarrow\mathscr{S}_{\phi}^1[1]=\mathcal{O}_{\phi}+\bar{\mathcal{O}}:(a,b)\mapsto a-b~,
\]
induces an isomorphism $\mathscr{H}^k(\mathscr{S}_{\phi}^\bullet[1])\cong\mathscr{H}^k(\mathscr{B}_{\phi}^\bullet),~\forall k\geq 0$.
\end{proposition}

\begin{proof}Note that $\mathscr{H}^{1}(\mathscr{B}_{\phi}^\bullet)=\mathbb{C}\oplus \mathbb{C}/\mathbb{C}(1,-1)\longrightarrow \mathbb{C}=\mathscr{H}^{1}(\mathscr{S}_{\phi}^\bullet[1]):(a,b)\mapsto a-b$ is an isomorphism.
\end{proof}

It follows that
\begin{equation}\label{eq-BC-hyper-Bphi}
H_{BC\phi}^{p,q}(X)\cong\mathbb{H}^{p+q}(M,\mathscr{L}^\bullet_{\phi}[1])\cong\mathbb{H}^{p+q}(M,\mathscr{S}^\bullet_{\phi}[1])\cong\mathbb{H}^{p+q}(M,\mathscr{B}_{\phi}^\bullet)~.
\end{equation}
Note that \eqref{eq-BC-hyper-Bphi} and Proposition \ref{prop-quasiiso-B-S-phi} is just a slight generalization of the result obtained by Schweitzer. In fact, Proposition \ref{prop-quasiiso-B-S-phi} reduce to \cite[Prop.\,4.3]{Sch07} when $\phi=0$.

Similarly, for fixed $p\geq 1$ and $q\geq 1$, we define a sheaf complex $\mathscr{L}_{\bar{\phi}}^\bullet$ as follows:
\begin{equation}
\left\{
\begin{array}{ll}
\mathscr{L}_{\bar{\phi}}^k= \bigoplus_{r+s=k,r<p,s<q}A^{r,s},~~~\text{for}~~~k\leq p+q-2, & \\
\mathscr{L}_{\bar{\phi}}^{k-1}= \bigoplus_{r+s=k,r\geq p,s\geq q}A^{r,s},~~~\text{for}~~~k\geq p+q.  & \\
\end{array} \right.
\end{equation}
The differential is given by
\begin{equation*}
\xymatrix@C=0.5cm{
  0 \ar[r] & \mathscr{L}_{\bar{\phi}}^0 \ar[rr]^{\Pi_{\mathscr{L}_{\bar{\phi}}^1}d_{\bar{\phi}}} && \mathscr{L}_{\bar{\phi}}^1 \ar[rr]^{\Pi_{\mathscr{L}_{\bar{\phi}}^2}d_{\bar{\phi}}} && \mathscr{L}_{\bar{\phi}}^2 \ar[r] & \cdots \\  \ar[r] & \mathscr{L}_{\bar{\phi}}^{p+q-3} \ar[rr]^{\Pi_{\mathscr{L}_{\bar{\phi}}^{p+q-2}}d_{\bar{\phi}}} && \mathscr{L}_{\bar{\phi}}^{p+q-2} \ar[rr]^{\partial_{\phi}\bar{\partial}} && \mathscr{L}_{\bar{\phi}}^{p+q-1} \ar[rr]^{d_{\bar{\phi}}} && \mathscr{L}_{\bar{\phi}}^{p+q} \ar[rr]^{d_{\bar{\phi}}} && \cdots .}
\end{equation*}
We have $\mathbb{H}^{p+q-1}(X,\mathscr{L}_{\bar{\phi}}^\bullet)\cong H^{p+q-1}(\mathscr{L}_{\bar{\phi}}^\bullet(X))=H_{BC\bar{\phi}}^{p,q}(X)$. The sheaf complex $\mathscr{L}_{\bar{\phi}}^\bullet$ has the following subcomplex
\begin{equation*}
\xymatrix@C=0.5cm{
\mathscr{S}_{\bar{\phi}}^\bullet:~\mathcal{O}\oplus\bar{\mathcal{O}}_{\phi} \ar[r]^{\quad\partial_{\phi}\oplus\bar{\partial}} & \Omega^{1}\oplus\bar{\Omega}^{1}_{\phi} \ar[r]^{\partial_{\phi}\oplus\bar{\partial}} & \Omega^{2}\oplus\bar{\Omega}^{2}_{\phi} \ar[r] & \cdots \\
 \ar[r] & \Omega^{p-1}\oplus\bar{\Omega}^{p-1}_{\phi}\ar[r]^{\qquad 0\oplus \bar{\partial}} & \bar{\Omega}^{p}_{\phi}\ar[r]^{\bar{\partial}} &\cdots  \ar[r]^{\bar{\partial}} & \bar{\Omega}^{q-1}_{\phi} \ar[r] & 0.}
\end{equation*}
\begin{proposition}\label{prop-quasiiso-S-L-phibar}
The inclusion $\mathscr{S}_{\bar{\phi}}^\bullet\hookrightarrow \mathscr{L}_{\bar{\phi}}^\bullet$ induces an isomorphism $\mathscr{H}^k(\mathscr{S}_{\bar{\phi}}^\bullet)\cong\mathscr{H}^k(\mathscr{L}_{\bar{\phi}}^\bullet),~\forall k\geq 0$, and we have
\begin{equation*}
\mathscr{H}^k(\mathscr{S}_{\bar{\phi}}^\bullet)\cong\mathscr{H}^k(\mathscr{L}_{\bar{\phi}}^\bullet)=\left\{
\begin{array}{ll}
\mathbb{C},~~~\text{for}~~~k=0,~p>1,~q>1, & \\
\mathcal{O},~~~\text{for}~~~k=0,~p=1,~q>1, & \\
\bar{\mathcal{O}}_{\phi},~~~\text{for}~~~k=0,~p>1,~q=1, & \\
\mathcal{O}\oplus\bar{\mathcal{O}}_{\phi},~~~\text{for}~~~k=0,~p=1,~q=1, & \\
\Omega^{p-1}/\partial_{\phi}\Omega^{p-2},~~~\text{for}~~~0<k=p-1\neq q-1,  & \\
\bar{\Omega}^{q-1}_{\phi}/\bar{\partial}\bar{\Omega}^{q-2}_{\phi},~~~\text{for}~~~0<k=q-1\neq p-1,  & \\
\Omega^{p-1}/\partial_{\phi}\Omega^{p-2}\oplus \bar{\Omega}^{p-1}_{\phi}/\bar{\partial}\bar{\Omega}^{p-2}_{\phi},~~~\text{for}~~~0<k=p-1=q-1,  & \\
0,~~~~~~\text{otherwise}.
\end{array} \right.
\end{equation*}
\end{proposition}

Consider the complex $\mathscr{B}_{\bar{\phi}}^\bullet$ which is a modification of $\mathscr{S}_{\bar{\phi}}^\bullet$ given by
\begin{equation*}
\xymatrix@C=0.5cm{
\mathscr{B}_{\bar{\phi}}^\bullet:\mathbb{C} \ar[r] & \mathcal{O}\oplus\bar{\mathcal{O}}_{\phi} \ar[r]^{\partial_{\phi}\oplus\bar{\partial}} & \Omega^{1}\oplus\bar{\Omega}^{1}_{\phi} \ar[r]^{\partial_{\phi}\oplus\bar{\partial}} & \Omega^{2}\oplus\bar{\Omega}^{2}_{\phi} \ar[r] & \cdots \\
 \ar[r] & \Omega^{p-1}\oplus\bar{\Omega}^{p-1}_{\phi}\ar[r]^{\qquad 0\oplus \bar{\partial}} & \bar{\Omega}^{p}_{\phi}\ar[r]^{\bar{\partial}} &\cdots  \ar[r]^{\bar{\partial}} & \bar{\Omega}^{q-1}_{\phi} \ar[r] & 0.}
\end{equation*}
\begin{proposition}\label{prop-quasiiso-B-S-phibar}
The natural map from $\mathscr{B}_{\bar{\phi}}^\bullet$ to $\mathscr{S}_{\bar{\phi}}^\bullet[1]$, where
\[
\mathscr{B}_{\bar{\phi}}^1=\bar{\mathcal{O}}_{\phi}\oplus\mathcal{O}\longrightarrow\mathscr{S}_{\bar{\phi}}^1[1]=\bar{\mathcal{O}}_{\phi}+\mathcal{O}:(a,b)\mapsto a-b~,
\]
induces an isomorphism $\mathscr{H}^k(\mathscr{S}_{\bar{\phi}}^\bullet[1])\cong\mathscr{H}^k(\mathscr{B}_{\bar{\phi}}^\bullet),~\forall k\geq 0$.
\end{proposition}

It follows that
\begin{equation}\label{eq-BCphibar-hyper-Bphibar}
H_{BC\bar{\phi}}^{p,q}(X)\cong\mathbb{H}^{p+q}(M,\mathscr{L}^\bullet_{\bar{\phi}}[1])\cong\mathbb{H}^{p+q}(M,\mathscr{S}^\bullet_{\bar{\phi}}[1])\cong\mathbb{H}^{p+q}(M,\mathscr{B}_{\bar{\phi}}^\bullet)~.
\end{equation}
\begin{remark}There are natural isomorphisms
\[
H_{BC\phi}^{p,0}(X)\cong H_{BC}^{p,0}(X_t):\sigma\mapsto e^{i_{\phi}}\sigma,~~H_{BC\bar{\phi}}^{0,q}(X)\cong H_{BC}^{0,q}(X_t):\sigma\mapsto e^{i_{\bar{\phi}}}\sigma,
\]
and note also that
\[
H_{BC\phi}^{0,q}(X)= H_{BC}^{0,q}(X),~~H_{BC\bar{\phi}}^{p,0}(X)\cong H_{BC}^{p,0}(X).
\]
\end{remark}

\subsection{The Bott-Chern cohomology on $X_t$}
\label{Relations between}
Let $X$ be a complex manifold and $X_t$ a small deformation (of $X$) whose complex structure is represented by a Beltrami differential $\phi\in A^{0,1}(X, T_{X}^{1,0})$, then by \cite[Th.\,4.3]{Xia19dDol} or \cite[Prop.\,2.13]{RZ18} we know that there are isomorphism of sheaves
\[
e^{i_{\phi}}:\Omega^{p}_{\phi}\longrightarrow \Omega^{p}_{X_t},~p=0,1,2,\cdots,n,
\]
which give rise to the following commutative diagram
\[
\xymatrix{
\mathscr{C}^\bullet:\mathbb{C}  \ar[r]^{(+,-)}  & \mathcal{O}_{\phi}\oplus\bar{\mathcal{O}}_{\phi} \ar[r]^{\partial\oplus\bar{\partial}}
\ar[d]^{id} &\Omega^{1}_{\phi}\oplus\bar{\Omega}^{1}_{\phi} \ar[r] \ar[d]^{e^{i_{\phi}}\oplus e^{i_{\bar{\phi}}}} &  \cdots \\
\mathscr{B}_{X_t}^\bullet:\mathbb{C} \ar[r]^{(+,-)}    & \mathcal{O}_{X_t}\oplus\bar{\mathcal{O}}_{X_t} \ar[r]^{\partial_t\oplus\bar{\partial}_{t}} & \Omega^{1}_{X_t}\oplus\bar{\Omega}^{1}_{X_t}  \ar[r] &  \cdots }
\]
\[
\xymatrix{
\ar[r] & \Omega^{q-1}_{\phi}\oplus\bar{\Omega}^{q-1}_{\phi}\ar[r]^{\qquad\partial\oplus 0} \ar[d]^{e^{i_{\phi}}\oplus e^{i_{\bar{\phi}}}} & \Omega^{q}_{\phi}\ar[r] \ar[d]^{e^{i_{\phi}}} & \cdots \ar[r]^{\partial} & \Omega^{p-1}_{\phi}\ar[r] \ar[d]^{e^{i_{\phi}}}& 0 \\
\ar[r] & \Omega^{q-1}_{X_t}\oplus\bar{\Omega}^{q-1}_{X_t}\ar[r]^{\qquad\partial_t\oplus 0} & \Omega^{q}_{X_t}\ar[r] &  \cdots \ar[r]^{\partial_t} & \Omega^{p-1}_{X_t}\ar[r] & 0 . }
\]
We see that
\[
H_{BC}^{p,q}(X_t)= H^{p+q-1}(\mathscr{L}_{X_t}^\bullet(M))\cong \mathbb{H}^{p+q-1}(M,\mathscr{L}_{X_t}^\bullet)\cong \mathbb{H}^{p+q}(M,\mathscr{B}_{X_t}^\bullet)\cong \mathbb{H}^{p+q}(M,\mathscr{C}^\bullet),
\]
where $M$ is the underlying smooth manifold of $X$ and $X_t$.

\subsection{The case of Aeppli cohomology}
The \emph{deformed Aeppli cohomology} can be defined as follows:
\begin{equation}\label{}
H^{p,q}_{A\phi}(X) := \frac{\ker \partial\bar{\partial}_{\phi} \cap  A^{p,q}(X)}{\Image d_{\phi} \cap A^{p,q}(X)},~~~~~~
H^{p,q}_{A\bar{\phi}}(X) := \frac{\ker \partial_{\phi}\bar{\partial} \cap A^{p,q}(X)}{\Image d_{\bar{\phi}} \cap A^{p,q}(X)},~~\forall p,q\geq 0~,
\end{equation}
and $h^{p,q}_{A\phi}:=\dim H^{p,q}_{A\phi}(X), h^{p,q}_{A\bar{\phi}}:=\dim H^{p,q}_{A\bar{\phi}}(X)$. The conjugation gives a natural isomorphism between $H^{p,q}_{A\phi}(X)$ and $H^{q,p}_{A\bar{\phi}}(X)$, we thus have $h^{p,q}_{A\phi}=h^{q,p}_{A\bar{\phi}}$.

For fixed $p\geq 0$ and $q\geq 0$, similar to the constructions for the Bott-Chern cohomology we define a sheaf complex which still denoted by $\mathscr{L}_{\phi}^\bullet$ as follows:
\begin{equation}
\left\{
\begin{array}{ll}
\mathscr{L}_{\phi}^k= \bigoplus_{r+s=k,r<p+1,s<q+1}A^{r,s},~~~\text{for}~~~k\leq p+q, & \\
\mathscr{L}_{\phi}^{k-1}= \bigoplus_{r+s=k,r\geq p+1,s\geq q+1}A^{r,s},~~~\text{for}~~~k\geq p+q+2.  & \\
\end{array} \right.
\end{equation}
The differential is given by
\begin{equation*}
\xymatrix@C=0.5cm{
  0 \ar[r] & \mathscr{L}_{\phi}^0 \ar[rr]^{\Pi_{\mathscr{L}_{\phi}^1}d_{\phi}} && \mathscr{L}_{\phi}^1 \ar[rr]^{\Pi_{\mathscr{L}_{\phi}^2}d_{\phi}} && \mathscr{L}_{\phi}^2 \ar[r] & \cdots \\  \ar[r] & \mathscr{L}_{\phi}^{p+q-1} \ar[rr]^{\Pi_{\mathscr{L}_{\phi}^{p+q}}d_{\phi}} && \mathscr{L}_{\phi}^{p+q} \ar[rr]^{\partial\bar{\partial}_{\phi}} && \mathscr{L}_{\phi}^{p+q+1} \ar[rr]^{d_{\phi}} && \mathscr{L}_{\phi}^{p+q+2} \ar[rr]^{d_{\phi}} && \cdots ,}
\end{equation*}
In particular, we find that
\begin{equation*}
\xymatrix@C=0.5cm{
\mathscr{L}_{\phi}^{p+q-1}=A^{p,q-1}\oplus A^{p-1,q} \ar[rr]^{\qquad\Pi_{\mathscr{L}_{\phi}^{p+q}}d_{\phi}} && \mathscr{L}_{\phi}^{p+q}=A^{p,q} \ar[rr]^{\partial\bar{\partial}_{\phi}\qquad} && \mathscr{L}_{\phi}^{p+q+1}=A^{p+1,q+1},}
\end{equation*}
and so $\mathbb{H}^{p+q}(X,\mathscr{L}_{\phi}^\bullet)\cong H^{p+q}(\mathscr{L}_{\phi}^\bullet(X))=H_{A\phi}^{p,q}(X)$. The other hypercohomology interpretations of the deformed Bott-Chern cohomology holds similarly for the deformed Aeppli cohomology.
The Hodge star operator induces the following duality between the deformed Bott-Chern cohomology and the deformed Aeppli cohomology~\cite[pp.\,10]{Sch07}:
\begin{equation}
H_{BC\phi}^{p,q}(X)\cong H_{A\phi}^{n-q,n-p}(X),~~\text{and}~~H_{BC\bar{\phi}}^{p,q}(X)\cong H_{A\bar{\phi}}^{n-q,n-p}(X)~.
\end{equation}

\section{Deformations of Bott-Chern classes}\label{Deformations of BC c}
Let $\pi: (\mathcal{X}, X) \to (B, 0)$ be a small deformation of a compact complex manifold $X$ such that for each $t\in B$ the complex structure on $X_t$ is represented by Beltrami differential $\phi(t)$. In this section, power series will always be written in homogenous form, e.g. we write $\phi(t)=\sum_k\phi_k$ where each $\phi_k$ is a homogeneous polynomial of degree $k$ with coefficients in $A^{0,1}(X,T^{1,0})$. The Bott-Chern Laplacian operator is defined as
\begin{equation}
\square_{BC}:=(\partial\bar{\partial})(\partial\bar{\partial})^*+(\partial\bar{\partial})^*(\partial\bar{\partial})+ (\bar{\partial}^*\partial)(\bar{\partial}^*\partial)^*+(\bar{\partial}^*\partial)^*(\bar{\partial}^*\partial)+ \bar{\partial}^*\bar{\partial} + \partial^*\partial,
\end{equation}
and the deformed Bott-Chern Laplacian operator is defined as
\begin{equation}
\square_{BC\phi}:=(\partial\bar{\partial}_{\phi})(\partial\bar{\partial}_{\phi})^*+(\partial\bar{\partial}_{\phi})^*(\partial\bar{\partial}_{\phi})+ (\bar{\partial}_{\phi}^*\partial)(\bar{\partial}_{\phi}^*\partial)^*+(\bar{\partial}_{\phi}^*\partial)^*(\bar{\partial}_{\phi}^*\partial)+ \bar{\partial}_{\phi}^*\bar{\partial}_{\phi} + \partial^*\partial,
\end{equation}
where $\phi=\phi(t)$.
Both $\square_{BC}$ and $\square_{BC\phi}$ are $4$-th order self-adjoint elliptic differential operator~\cite{Sch07,MK71}. We have
\begin{equation}
\mathcal{H}_{BC}:=\ker\square_{BC}= \ker\partial\cap\ker\bar{\partial}\cap\ker (\partial\bar{\partial})^*
\end{equation}
and the following orthogonal direct sum decomposition holds:
\begin{equation}\label{eq-BC-Hodgedecomposition}
A^{\bullet,\bullet}(X)=\ker\square_{BC}\oplus \Image\partial\bar{\partial}\oplus (\Image\partial^*+\Image\bar{\partial}^*),
\end{equation}
which is equivalent to the existence of the Green operator $G_{BC}$ such that
\[
1=\mathcal{H}_{BC}+\square_{BC}G_{BC}.
\]
The same is true for the deformed Bott-Chern Laplacian operator $\square_{BC\phi}$.
It follows from \eqref{eq-BC-Hodgedecomposition} that
\begin{equation}\label{eq-ker-ddbar*}
\ker (\partial\bar{\partial})^*=\mathcal{H}_{BC}\oplus(\Image \partial^*+\Image\bar{\partial}^*)~\text{and}~\ker d=\mathcal{H}_{BC}\oplus \Image\partial\bar{\partial}~.
\end{equation}
The Aeppli Laplacian operator is defined as
\begin{equation}
\square_{A}:=(\partial\bar{\partial})^*(\partial\bar{\partial})+(\partial\bar{\partial})(\partial\bar{\partial})^*+ (\bar{\partial}\partial^*)^*(\bar{\partial}^*\partial)+(\bar{\partial}\partial^*)(\bar{\partial}\partial^*)^*+ \bar{\partial}\bar{\partial}^* + \partial\partial^*,
\end{equation}
and we have correspondingly
\begin{equation}\label{eq-A-Hodgedecomposition}
A^{\bullet,\bullet}(X)=\ker\square_{A}\oplus \Image(\partial\bar{\partial})^*\oplus (\Image\partial+\Image\bar{\partial}),
\end{equation}
or $1=\mathcal{H}_{A}+\square_{A}G_{A}$ where $G_{A}$ is the Green operator for $\square_{A}$. Since for any $x\in A^{\bullet,\bullet}(X)$,
we have $\square_{BC}G_{BC}\partial\bar{\partial}x=\partial\bar{\partial}x$ and $\square_{BC}\partial\bar{\partial}G_{A}x=\partial\bar{\partial}x$ which implies
\begin{equation}\label{eq-GBC-ddbar}
G_{BC}\partial\bar{\partial}=\partial\bar{\partial}G_{A}~.
\end{equation}
Similarly, we have
\begin{equation}\label{eq-GA-ddbar*}
(\partial\bar{\partial})^*G_{BC}=G_{A}(\partial\bar{\partial})^*~.
\end{equation}

Let $\varphi\in A^{p,q}(X)$ and $G_{BC}: A^{p,q}(X)\to A^{p,q}(X)$ be the Green operator, then for $k\geq 2$ we have
\begin{equation}
\|G_{BC}\varphi\|_{k+\alpha}\leq C \|\varphi\|_{k-4+\alpha} ,
\end{equation}
where $C>0$ is independent of $\varphi$ and $\|\cdot\|_{k+\alpha}$ is the H\"older norm.

We have the following observation:
\begin{proposition}\label{prop-unique small solution}
$1$. $\forall \sigma \in A^{p,q}(X)$, if $d_{\phi(t)} \sigma = d\sigma - \mathcal{L}_{\phi(t)}^{1,0} \sigma =0$ and $(\partial\bar{\partial})^*\sigma=0$, then we must have
\[
\sigma = \mathcal{H}_{BC}\sigma - G_{BC}A \partial i_{\phi(t)} \sigma ,
\]
where $\mathcal{H}_{BC}: A^{p,q}(X)\to \mathcal{H}^{p,q}_{BC}(X)$ is the projection operator to harmonic space and $A:=\bar{\partial}^*\partial\partial^*+\bar{\partial}^*$.\\
$2$. For any fixed $\sigma_0\in \mathcal{H}^{p,q}_{BC}(X)$, the equation
\begin{equation}\label{Kuranishi eq}
\sigma=\sigma_0 - G_{BC}A \partial i_{\phi(t)} \sigma ,
\end{equation}
has an unique solution given by $\sigma=\sum_{k} \sigma_k\in A^{p,q}(X)$ and $\sigma_{k}=-G_{BC}A\sum_{i+j=k} \partial i_{\phi_j} \sigma_i $ for $|t|$ small where each $\sigma_k$ is a homogeneous polynomial of degree $k$ with coefficients in $A^{p,q}(X)$.
\end{proposition}

\begin{proof}The first assertion follows from the Hodge decomposition:
\[
\sigma = \mathcal{H}_{BC}\sigma + G_{BC}\square_{BC} \sigma
       = \mathcal{H}_{BC}\sigma + G_{BC}A \mathcal{L}_{\phi(t)}^{1,0} \sigma
       = \mathcal{H}_{BC}\sigma - G_{BC}A \partial i_{\phi(t)} \sigma~,
\]
where we have used the fact that $d_{\phi(t)} \sigma =0 \Leftrightarrow  \partial\sigma= \bar{\partial}\sigma + \partial i_{\phi(t)} \sigma =0$.

For the second assertion, substitute $\sigma=\sigma(t)=\sum_{k}$ in \eqref{Kuranishi eq}, we have
\begin{equation}\label{formal solutions}
\left\{
\begin{array}{ll}
\sigma_1 &= - G_{BC}A \partial i_{\phi_1} \sigma_0,  \\
\sigma_2 &= - G_{BC}A (\partial i_{\phi_2} \sigma_0 + \partial i_{\phi_1} \sigma_1 ),\\
         &\cdots,  \\
\sigma_k &= - G_{BC}A\sum_{i+j=k} \partial i_{\phi_j} \sigma_i ,~~\forall k > 0.  \\
\end{array} \right.
\end{equation}
For the convergence of $\sigma(t)$, we note that
\begin{equation}
\|\sigma_{j}\|_{k+\alpha}=\| G_{BC}A\sum_{a+b=j} \partial i_{\phi_a} \sigma_b \|_{k+\alpha} \leq C \sum_{a+b=j} \|\phi_{a}\|_{k+\alpha} \|\sigma_{b}\|_{k+\alpha},
\end{equation}
for some constant $C$ depends only on $k$ and $\alpha$.
Now it is left to show the uniqueness. Let $\sigma$ and $\sigma'$ be two solutions to $\sigma=\sigma_0 - G_{BC}A \partial i_{\phi} \sigma$ and set $\tau=\sigma-\sigma'$. Then $\tau=-G_{BC}A \partial i_{\phi} \tau$, we have
\begin{equation}
\|\tau\|_{k+\alpha}\leq c \|\phi (t)\|_{k+\alpha}  \|\tau\|_{k+\alpha},
\end{equation}
for some constant $c>0$. When $|t|$ is sufficiently small, $\|\phi (t)\|_{k+\alpha}$ is also small. Hence we must have $\tau=0$.
For smoothness of the solution, note that we have
\[
\square_{BC}\sigma=-\square_{BC}G_{BC}A \partial i_{\phi(t)} \sigma=-(1-\mathcal{H}_{BC\phi(t)} )A \partial i_{\phi(t)} \sigma=-A \partial i_{\phi(t)} \sigma,
\]
which implies
\begin{equation}\label{regularity eq}
\square_{BC}\sigma + (\bar{\partial}^*\partial\partial^*+\bar{\partial}^*)\partial i_{\phi(t)} \sigma= 0,
\end{equation}
which is a standard elliptic equation for small $t$.
\end{proof}
Note that the solution $\sigma$ of \eqref{Kuranishi eq} automatically satisfies $(\partial\bar{\partial})^*\sigma=0$ in view of \eqref{eq-GA-ddbar*}.

\begin{lemma}\label{lem-ddbar*-closed-rep of deformed-BC-coho}
The natural map
\begin{equation}\label{eq-natural map}
\frac{\ker(\partial\bar{\partial})^*\cap\ker d_{\phi(t)} \cap A^{p,q}(X)}{\ker(\partial\bar{\partial})^*\cap\Image\partial\bar{\partial}_{\phi(t)} \cap A^{p,q}(X)}\longrightarrow H_{BC\phi(t)}^{p,q}(X)
\end{equation}
is an isomorphism.
\end{lemma}

\begin{proof}By (the deformed version of) \eqref{eq-BC-Hodgedecomposition} and \eqref{eq-ker-ddbar*}, we have the following orthogonal direct sum decomposition
\[
A^{p,q}(X)=\left(\ker d_{\phi(t)}\cap\ker(\partial\bar{\partial})^*\right)\oplus
(\Image\partial^*+\Image\bar{\partial}_{\phi(t)}^*+\Image\partial\bar{\partial})
\]
which implies
\[
\ker d_{\phi(t)}\cap A^{p,q}(X)=\left(\ker d_{\phi(t)}\cap\ker(\partial\bar{\partial})^*\right)\oplus
\left(\ker d_{\phi(t)}\cap(\Image\partial^*+\Image\bar{\partial}_{\phi(t)}^*+\Image\partial\bar{\partial})\right),
\]
and
\[
\Image\partial\bar{\partial}_{\phi(t)}\cap A^{p,q}(X)=\left(\Image\partial\bar{\partial}_{\phi(t)}\cap\ker(\partial\bar{\partial})^*\right)\oplus
\left(\Image\partial\bar{\partial}_{\phi(t)}\cap(\ker(\partial\bar{\partial}_{\phi(t)})^*+\Image\partial\bar{\partial})\right).
\]
Moreover, for any $x\in \Image\partial\bar{\partial}$, there exists unique $y\in\ker d_{\phi(t)}$ and unique $z\in (\Image\partial^*+\Image\bar{\partial}_{\phi(t)}^*)$ such that $x=y+z$. This defines a surjective homomorphism
\[
\Image\bar{\partial}\bar{\partial}\longrightarrow\ker d_{\phi(t)}\cap(\Image\partial^*+\Image\bar{\partial}_{\phi(t)}^*+\Image\partial\bar{\partial}):
x\longmapsto y,
\]
with kernel equal to $\Image\partial\bar{\partial}\cap(\Image\partial^*+\Image\bar{\partial}_{\phi(t)}^*)$. It follows that
\[
\ker d_{\phi(t)}\cap(\Image\partial^*+\Image\bar{\partial}_{\phi(t)}^*+\Image\partial\bar{\partial})
\cong\frac{\Image\partial\bar{\partial}}{\Image\partial\bar{\partial}\cap(\Image\partial^*+\Image\bar{\partial}_{\phi(t)}^*)}.
\]
Similarly, we have
\[
\Image\partial\bar{\partial}_{\phi(t)}\cap(\ker(\partial\bar{\partial}_{\phi(t)})^*+\Image\partial\bar{\partial})
\cong\frac{\Image\partial\bar{\partial}}{\Image\partial\bar{\partial}\cap\ker(\partial\bar{\partial}_{\phi(t)})^*}.
\]
Hence,
\[
H_{BC\phi(t)}^{p,q}(X)\cong
\frac{\ker(\partial\bar{\partial})^*\cap\ker d_{\phi(t)}}{\ker(\partial\bar{\partial})^*\cap\Image\partial\bar{\partial}_{\phi(t)} } \oplus
\frac{\Image\partial\bar{\partial}\cap\ker(\partial\bar{\partial}_{\phi(t)})^*}{\Image\partial\bar{\partial}\cap(\Image\partial^*+\Image\bar{\partial}_{\phi(t)}^*)}
.\]
We claim $\Image\partial\bar{\partial}\cap\ker(\partial\bar{\partial}_{\phi(t)})^*=0$. Indeed, let $\sigma\in\ker\partial\bar{\partial}\cap\ker(\partial\bar{\partial}_{\phi(t)})^*$, then it follows from the same proof of Proposition \ref{prop-unique small solution} that $\sigma$ is the solution of the equation
\[
\sigma=\sigma_0+G_{BC}\partial\bar{\partial}(\partial\mathcal{L}_{\phi(t)}^{1,0})^*\sigma,\quad \sigma_0:=\mathcal{H}_{BC}\sigma
\]
which is uniquely determined by $\sigma_0$. If $\sigma\in \Image\partial\bar{\partial}\cap\ker(\partial\bar{\partial}_{\phi(t)})^*$, then $\sigma_0=\mathcal{H}_{BC}\sigma=0\Rightarrow\sigma=0$.
\end{proof}

\begin{proposition}\label{prop-kerddbar*-imageddbarphi}
$1$. For any fixed $t\in B$, the following homomorphism
\[
g_t:\ker\partial\bar{\partial}\cap A^{p,q}(X)\longrightarrow \ker(\partial\bar{\partial})^*\cap\Image\partial\bar{\partial}_{\phi(t)}\cap A^{p+1,q+1}(X):x_0\longmapsto \partial\bar{\partial}_{\phi(t)}x(t),
\]
is surjective with $\ker g_t=\ker\partial\bar{\partial}\cap\left(\ker\partial\bar{\partial}_{\phi(t)}+\Image(\partial\bar{\partial})^*\right)\cap A^{p,q}(X)$, where $x(t)$ is the unique solution of $x(t)=x_0 + (\partial\bar{\partial})^*G_{BC} \partial i_{\phi(t)} \partial x(t)$.\\
$2$. Let $\hat{g}_t:\mathcal{H}^{p,q}_{BC}(X)\longrightarrow \ker(\partial\bar{\partial})^*\cap\Image\partial\bar{\partial}_{\phi(t)}\cap A^{p+1,q+1}(X)$ be the restriction of $g_t$ on $\mathcal{H}^{p,q}_{BC}(X)$, then $\hat{g}_t$ is surjective with
$\ker \hat{g}_t=\mathcal{H}^{p,q}_{BC}(X)\cap\left(\ker\partial\bar{\partial}_{\phi(t)}+\Image(\partial\bar{\partial})^*\right)$. Moreover, we have
\begin{align*}
\dim\mathcal{H}^{p,q}_{BC}(X)=&\dim\ker\partial\bar{\partial}_{\phi(t)}\cap(\mathcal{H}^{p,q}_{BC}(X)+\Image(\partial\bar{\partial})^*)\cap A^{p,q}(X)\\+
&\dim \ker(\partial\bar{\partial})^*\cap\Image\partial\bar{\partial}_{\phi(t)}\cap A^{p+1,q+1}(X).
\end{align*}
\end{proposition}
\begin{proof}$1$. Let $x\in A^{p,q}(X)$, then by Hodge decomposition we have
\[
\partial\bar{\partial}_{\phi(t)}x=\partial\bar{\partial}x-\partial i_{\phi(t)}\partial x
=\partial\bar{\partial}x-\mathcal{H}_{BC}\partial i_{\phi(t)}\partial x-\square_{BC}G_{BC}\partial i_{\phi(t)}\partial x,
\]
thus
\[
\partial\bar{\partial}_{\phi(t)}x\in\ker(\partial\bar{\partial})^*\Leftrightarrow \partial\bar{\partial}x-\partial\bar{\partial}(\partial\bar{\partial})^*G_{BC}\partial i_{\phi(t)}\partial x=0.
\]
Set $x_0=x-(\partial\bar{\partial})^*G_{BC}\partial i_{\phi(t)}\partial x$, then $x$ is a solution to the equation $x=x_0 + (\partial\bar{\partial})^*G_{BC}\partial i_{\phi(t)}\partial x$ which is uniquely determined by $x_0$ in view of the proof of Proposition \ref{prop-unique small solution}.

It is left to show $\ker g_t=\ker\partial\bar{\partial}\cap\left(\ker\partial\bar{\partial}_{\phi(t)}+\Image(\partial\bar{\partial})^*\right)\cap A^{p,q}(X)$. In fact, obviously we have
$\ker g_t\subseteq\ker\partial\bar{\partial}\cap\left(\ker\partial\bar{\partial}_{\phi(t)}+\Image(\partial\bar{\partial})^*\right)\cap A^{p,q}(X)$. Conversely, let us consider the following surjective homomorphism
\begin{align*}
\ker&\partial\bar{\partial}_{\phi(t)}\longrightarrow \ker\partial\bar{\partial}\cap\left(\ker\partial\bar{\partial}_{\phi(t)}+\Image(\partial\bar{\partial})^*\right)\\
&x\longmapsto x-(\partial\bar{\partial})^*\partial\bar{\partial}G_Ax
=x-(\partial\bar{\partial})^*G_{BC}\partial\bar{\partial}x=x-(\partial\bar{\partial})^*G_{BC}\partial i_{\phi(t)}\partial x,
\end{align*}
whose kernel is $\ker\partial\bar{\partial}_{\phi(t)}\cap\Image(\partial\bar{\partial})^*=0$ by Proposition \ref{prop-unique small solution}. Its inverse is given by
\[
\ker\partial\bar{\partial}\cap\left(\ker\partial\bar{\partial}_{\phi(t)}+\Image(\partial\bar{\partial})^*\right)
\longrightarrow \ker\partial\bar{\partial}_{\phi(t)}:
x_0\longmapsto x(t),
\]
where $x(t)$ is the unique solution of $x(t)=x_0 + (\partial\bar{\partial})^*G_{BC}\partial i_{\phi(t)}\partial x(t)$.
So let $x_0\in\ker\partial\bar{\partial}\cap\left(\ker\partial\bar{\partial}_{\phi(t)}+\Image(\partial\bar{\partial})^*\right)$ and $x(t)$ be the unique solution of $x(t)=x_0 + (\partial\bar{\partial})^*G_{BC}\partial i_{\phi(t)}\partial x(t)$, we must have $x(t)\in \ker\partial\bar{\partial}_{\phi(t)}\Rightarrow x_0\in\ker g_t$. \\

$2$. It can be proved in the same way that $\ker \hat{g}_t=\mathcal{H}^{p,q}_{BC}(X)\cap\left(\ker\partial\bar{\partial}_{\phi(t)}+\Image(\partial\bar{\partial})^*\right)$. To show $\hat{g}_t$ is surjective it is enough to show
\[
\frac{\mathcal{H}^{p,q}_{BC}(X)}{\mathcal{H}^{p,q}_{BC}(X)\cap\left(\ker\partial\bar{\partial}_{\phi(t)}+\Image(\partial\bar{\partial})^*\right)}\cong
\frac{\ker\partial\bar{\partial}\cap A^{p,q}(X)}{\ker\partial\bar{\partial}\cap\left(\ker\partial\bar{\partial}_{\phi(t)}+\Image(\partial\bar{\partial})^*\right)\cap A^{p,q}(X)}.
\]
Indeed, we have
\[
\ker\partial\bar{\partial}=\mathcal{H}^{p,q}_{BC}(X)\oplus\left\{\ker\partial\bar{\partial}\cap\left(\Image(\partial\bar{\partial})+ \Image\partial^*+ \Image\bar{\partial}^*\right)\right\},
\]
and
\begin{align*}
&\ker\partial\bar{\partial}\cap\left(\ker\partial\bar{\partial}_{\phi(t)}+\Image(\partial\bar{\partial})^*\right)\cong\ker\partial\bar{\partial}_{\phi(t)}\\
=&\left\{\ker\partial\bar{\partial}_{\phi(t)}\cap (\mathcal{H}^{p,q}_{BC}(X)+\Image(\partial\bar{\partial})^*)\right\}\oplus\\
&\left\{\ker\partial\bar{\partial}_{\phi(t)}\cap\left[\Image(\partial\bar{\partial}_{\phi(t)})^*+ \ker\partial\bar{\partial}\cap(\Image\partial\bar{\partial} +\Image\partial^*+\Image\bar{\partial}^*)\right]\right\}\\
\cong& \left\{\mathcal{H}^{p,q}_{BC}(X)\cap\left(\ker\partial\bar{\partial}_{\phi(t)}+\Image(\partial\bar{\partial})^*\right)\right\}
\oplus \left\{\ker\partial\bar{\partial}\cap\left(\Image(\partial\bar{\partial})+ \Image\partial^*+ \Image\bar{\partial}^*\right)\right\}.
\end{align*}
\end{proof}
\begin{remark}It can be proved in a similar way that $g_t$ when restricted on $\mathcal{H}^{p,q}_{A}(X)$ is also surjective with kernel equal to
$\mathcal{H}^{p,q}_{A}(X)\cap\left(\ker\partial\bar{\partial}_{\phi(t)}+\Image(\partial\bar{\partial})^*\right)$ and
\begin{align*}
\dim\mathcal{H}^{p,q}_{A}(X)=&\dim\ker\partial\bar{\partial}_{\phi(t)}\cap(\mathcal{H}^{p,q}_{A}(X)+\Image(\partial\bar{\partial})^*)\cap A^{p,q}(X)\\+
&\dim \ker(\partial\bar{\partial})^*\cap\Image\partial\bar{\partial}_{\phi(t)}\cap A^{p+1,q+1}(X).
\end{align*}
\end{remark}

\begin{definition}\label{def-V_t}
For any $t\in B$ and a vector subspace $V=\mathbb{C}\{ \sigma_0^1, \cdots, \sigma_0^N \}\subseteq \mathcal{H}_{BC}^{p,q}(X)$, we set
\begin{align*}\label{V_E,t}
V_{t}:=
&\{ \sum_{l=1}^N a_l\sigma_0^l\in V \mid  (a_1, \cdots, a_N)\in \mathbb{C}^{N}~\text{s.t.}~\sigma(t)\in\ker d_{\phi(t)},\\
&\text{where}~\sigma (t)=\sum_{k} \sigma_k~\text{with}~ \sigma_0=\sum_l a_l\sigma_0^l~\text{and}~ \sigma_{k}=-G_{BC}A\sum_{i+j=k} \partial i_{\phi_j} \sigma_i,~\forall k\neq 0 \}.
\end{align*}
Note that $V_{t}$ consists of those vectors of the form $\sum_l a_l\sigma_0^l$ such that the coefficients $a_l$ satisfy the following linear equation:
\[
\sum_{l=1}^N a_ld_{\phi(t)}\sigma^l(t)=0,
\]
where $\sigma^l(t)= \sum_{k} \sigma_k^l$ with $\sigma_{k}^l=-G_{BC}A\sum_{i+j=k} \partial i_{\phi_j}\sigma_i,~\forall k\neq 0$.
\end{definition}

\begin{definition}\label{def-f_t-g_t}
We set
\begin{align*}
f_t: &V_{t} \longrightarrow \frac{\ker(\partial\bar{\partial})^*\cap\ker d_{\phi(t)} \cap A^{p,q}(X)}{\ker(\partial\bar{\partial})^*\cap\Image\partial\bar{\partial}_{\phi(t)} \cap A^{p,q}(X)}\cong H_{BC\phi(t)}^{p,q}(X),\\
&\sigma_0\longmapsto \sigma(t)=\sum_{k} \sigma_k,~\text{where}~\sigma_{k}=-G_{BC}A\sum_{i+j=k} \partial i_{\phi_j} \sigma_i,~\forall k\neq 0.
\end{align*}
\end{definition}

\begin{proposition}\label{prop V-E-t}
If $V= \mathcal{H}_{BC}^{p,q}(X)$, then $f_t$ is surjective.
\end{proposition}

\begin{proof}By Proposition \ref{prop-unique small solution}, The map
\begin{align*}
\tilde{f}_t: &V_{t} \longrightarrow \ker(\partial\bar{\partial})^*\cap\ker d_{\phi(t)} \cap A^{p,q}(X),\\
&\sigma_0\longmapsto \sigma(t)=\sum_{k} \sigma_k,~\text{where}~\sigma_{k}=-G_{BC}A\sum_{i+j=k} \partial i_{\phi_j} \sigma_i,~\forall k\neq 0,
\end{align*}
is an isomorphism.
\end{proof}

\begin{theorem}\label{Deformation of B-C cohomology classes}
Let $X$ be a compact complex manifold and $\pi: (\mathcal{X}, X) \to (B, 0)$ a small deformation of $X$ such that for each $t\in B$ the complex structure on $X_t$ is represented by Beltrami differential $\phi(t)$. For any $p,q\geq 0$, let $V=\mathbb{C}\{ \sigma_0^1, \cdots, \sigma_0^N \}$ be a linear subspace of $\mathcal{H}_{BC}^{p,q}(X)$ and $\sigma^l(t)=\tilde{f}_t\sigma_{0}^l,~l=1, \cdots, N$. Define a subset $B(V)$ of $B$ by
\[
B(V):=\{t\in \mathcal{B}\mid d_{\phi(t)}\sigma^l(t)=0, l=1, \cdots, N\},
\]
Then $B(V)$ are analytic subsets of $B$ and we have
\begin{equation}\label{B(V)-discription}
B(V)=\{t\in B\mid \dim V = \dim \Image f_t +\dim \ker f_t\}.
\end{equation}
In particular, we have
\begin{equation}\label{B'-discription}
B'=B(\mathcal{H}_{BC}^{p,q}(X))=\{t\in B\mid \dim H_{BC}^{p,q}(X)=\dim H_{BC\phi(t)}^{p,q}(X)+\dim \ker f_t\}.
\end{equation}
\end{theorem}

\begin{proof}First, let $\{U_{\alpha}\} $ be a finite open cover of $X$ and $ u^\alpha_1, u^\alpha_{2}, \cdots, $ a local unitary frames of $p+q+1$-forms on the $U_{\alpha}$, then $\forall l=1, \cdots, N$, we have
\[
d_{\phi(t)}\sigma^l(t)=0 \Leftrightarrow a_{j}^{l\alpha}(t):=< d_{\phi(t)}\sigma^l(t)\mid_{U_{\alpha}} ,u^\alpha_j>=0,~\forall j,~\alpha,
\]
where $<\cdot,\cdot>$ is the $L^2$-inner product on the space $A^{p+q+1}(U_{\alpha})$. We see that each $a_{j}^{l\alpha}(t)$ is holomorphic in $t$ and so
\[
B(V)=\{t\in B\mid a_{j}^{l\alpha}(t)=0, \forall j, l, \alpha\}
\]
is an analytic subset of $B$.

Furthermore, note that
\[
 t\in B(V) \Leftrightarrow V_{t}=V .
\]
So \eqref{B(V)-discription} follows from the fact that $\dim V_{t}=\dim \Image f_t +\dim \ker f_t$. If $V= \mathcal{H}^{0,q}(X,E)$, then $f_t: V_{t} \to H^{0,q}_{\bar{\partial}_{\phi(t)}}(X,E)$ is surjective by Proposition \ref{prop V-E-t} and \eqref{B'-discription} follows.
\end{proof}
\begin{remark}\label{rk-analyticity-Vt}
From the above proof, we can see that $V_{t}\subseteq V$ varies holomorphically with $t$. In fact, in the notations of Definition \ref{def-V_t}, $V_{t}$ consists of those vectors of the form $\sum_l a_l\sigma_0^l$ such that the coefficients $a_l$ satisfy
\[
\sum_{l=1}^N a_l\cdot a_{j}^{l\alpha}(t)=0,\quad j=1,\cdots,m;~\forall \alpha~,
\]
where $a_{j}^{l\alpha}(t)$ are holomorphic functions in $t$. In particular, $\{t\in B\mid \dim V_t\geq k\}$ is an analytic subset of $B$ for any nonnegative integer $k$.

For the same reason, $\{t\in B\mid\dim\ker\partial\bar{\partial}_{\phi(t)}\cap(\mathcal{H}^{p,q}_{BC}(X)+\Image(\partial\bar{\partial})^*)\cap A^{p,q}(X)\geq k\}$ is also an analytic subset of $B$ for any nonnegative integer $k$. It follows from this and Proposition \ref{prop-kerddbar*-imageddbarphi} that $\{t\in B\mid\dim \ker(\partial\bar{\partial})^*\cap\Image\partial\bar{\partial}_{\phi(t)}\cap A^{p,q}(X)\leq k\}$ is an analytic subset of $B$ for any nonnegative integer $k$.
\end{remark}

\begin{definition}\label{def-deformation-Bott-Chern}
Let $\pi: (\mathcal{X}, X)\to (B,0)$ be a deformation of a compact complex manifold $X$ such that for each $t\in B$ the complex structure on $X_t$ is represented by Beltrami differential $\phi(t)$. Given $y\in \ker d\cap A^{p,q}(X)$ and $T\subseteq B$, which is a complex subspace of $B$ containing $0$, a (Bott-Chern) \emph{deformation} of $y$ (w.r.t. $\pi: (\mathcal{X}, X)\to (B,0)$ ) on $T$ is a family of $(p,q)$-forms $\sigma (t)$ such that
\begin{itemize}
  \item[1.] $\sigma (t)$ is holomorphic in $t\in T$ and $\sigma (0) = y$;
  \item[2.] $d_{\phi(t)}\sigma (t) = 0,~\forall t\in T$.
\end{itemize}
A deformation of $[y]\in \mathcal{H}_{BC}^{p,q}(X)$ (w.r.t. $\pi$) on $T$ is a triple $(y,\sigma (t),T)$ which consisting of a representative $y\in[y]$ and a deformation $\sigma (t)$ of $y$ (w.r.t. $\pi$) on $T$. Two deformations $(y,\sigma (t),T)$ and $(y',\sigma' (t),T)$ of $[y]$ on $T$ are \emph{equivalent} if
\[
[\sigma (t) - \sigma' (t)] = 0 \in H_{BC}^{p,q}(X),~\forall t\in T.
\]
A deformation $\sigma (t)$ of $y$ on $T$ is called \emph{canonical} if
\[
\sigma(t)=\sigma_0 - G_{BC}(\bar{\partial}^*\partial\partial^*+\bar{\partial}^*) \partial i_{\phi(t)} \sigma(t),\quad \forall t\in T .
\]
By Proposition \ref{prop-unique small solution}, canonical deformation is unique on its existence domain.

For a given small deformation $\pi: (\mathcal{X}, X)\to (B,0)$ with smooth $B$, we say $y\in \ker d\cap A^{p,q}(X)$ is \emph{(canonically) unobstructed w.r.t. $\pi$} if a (canonical) deformation of $y$ (w.r.t. $\pi$) exists on $B$ and a class $\alpha\in H_{BC}^{p,q}(X)$ is \emph{(canonically) unobstructed w.r.t. $\pi$} if there is a $y\in\alpha$ such that $y$ is (canonically) unobstructed w.r.t. $\pi$. If every Bott-Chern classes in $H_{BC}^{p,q}(X)$ have canonically unobstructed deformation w.r.t. $\pi$, then we say the \emph{deformations of classes in $H_{BC}^{p,q}(X)$ is canonically unobstructed w.r.t. $\pi$}. If these holds for any small deformation of $X$, we will drop the term ``w.r.t. $\pi$''. For example, we say $y\in \ker d\cap A^{p,q}(X)$ is \emph{(canonically) unobstructed} if for any small deformation $\pi: (\mathcal{X}, X)\to (B,0)$ with smooth $B$ there is a (canonical) deformation of $y$ on $B$.
\end{definition}

Although a Bott-Chern deformation $\sigma (t)$ of $y\in \ker d\cap A^{p,q}(X)$ can also be viewed as a Dolbeault deformation with the additional requirement $\partial\sigma (t)=0$, the ways we identify deformations in these two cases is very different.
We want to point out another difference between the deformation theory of Dolbeault cohomology~\cite{Xia19dDol} and that of Bott-Chern cohomology. Let $\sigma^{BC} (t)$ and $\sigma^{D} (t)$ be the canonical Bott-Chern/Dolbeault deformation of $y\in \ker d\cap A^{p,q}(X)$ respectively, it is known that $\mathcal{L}_{\phi(t)}^{1,0}\sigma^{D}(t)\in\ker\bar{\partial}$ for any $t\in B$, see \cite[Prop.\,5.2]{Xia19dDol}. This seems does not hold for the Bott-Chern deformation in general. More precisely, it is not guaranteed that $\mathcal{L}_{\phi(t)}^{1,0}\sigma^{BC} (t)\in\ker d_{\phi(t)}$.

In the remainder of this section, We confine ourselves to sketching the essential points of the deformation theory of Bott-Chern cohomology. Since this part of the theory is very similar to the case of Dolbeault cohomology, the proofs will be omitted.

A notable consequence of the deformation theory for Bott-Chern classes is the following
\begin{theorem}\label{th-analytic subset}
Let $\pi: (\mathcal{X}, X)\to (B,0)$ be a deformation of a compact complex manifold $X$ such that for each $t\in B$ the complex structure on $X_t$ is represented by Beltrami differential $\phi(t)$. Then the set $\{t\in B\mid \dim H_{BC\phi(t)}^{p,q}(X)\geq k\}$ is an analytic subset of $B$ for any nonnegative integer $k$.
\end{theorem}
\begin{proof}It follows from Proposition \ref{prop V-E-t} that
\begin{align*}
&\{t\in B\mid \dim H_{BC\phi(t)}^{p,q}(X)\geq k\}\\
=&\{t\in B\mid \dim V_t/\ker f_t\geq k\}\\
=&\{t\in B\mid \dim V_t-\dim\left(\ker\partial\bar{\partial}^*\cap\Image\partial\bar{\partial}_{\phi(t)}\right)\geq k\}.
\end{align*}
The conclusion then follows from Remark \ref{rk-analyticity-Vt}.
\end{proof}

The canonical deformations has the following properties:
\begin{theorem}\label{universal and uniqueness of the canonical deformation}
Let $\pi: (\mathcal{X}, X)\to (B,0)$ be a deformation of a compact complex manifold $X$ such that for each $t\in B$ the complex structure on $X_t$ is represented by Beltrami differential $\phi(t)$.
\begin{itemize}
	\item[$(i)$] Assume $S$ is an analytic subset of $B$ with $0\in S$ and $y\in \ker d\cap A^{p,q}(X)$. If the canonical deformation of $y$ exists on $S$ then we must have $S\subseteq B(\mathbb{C}\mathcal{H}_{BC}y)$;
	\item[$(ii)$] For any deformed Bott-Chern cohomology class $[u]\in H_{BC\phi(t)}^{p,q}(X)$, there exists $\sigma_0\in \mathcal{H}_{BC}^{p,q}(X)$ such that $[u]=[\sigma(t)]$ where $\sigma(t)$ is the canonical deformation of $\sigma_0$.
\end{itemize}
\end{theorem}
\begin{proof}$(i)$ follows from Theorem \ref{Deformation of B-C cohomology classes} and $(ii)$ follows from Lemma \ref{lem-ddbar*-closed-rep of deformed-BC-coho}.
\end{proof}

We end this section with the following result which is of particular interests.
\begin{theorem}\label{thm-2nd-main}
Let $\pi: (\mathcal{X}, X)\to (B,0)$ be a small deformation of the compact complex manifold $X$ such that for each $t\in B$ the complex structure on $X_t$ is represented by Beltrami differential $\phi(t)$. For each $(p,q)\in \mathbb{N}\times \mathbb{N}$, set
\[
v^{p,q}_t:=\dim H_{BC}^{p,q}(X)-\dim \ker d_{\phi(t)}\cap\ker(\partial\bar{\partial})^*\cap A^{p,q}(X)\geq 0,
\]
and
\[
u^{p,q}_t:=\dim H_{BC}^{p,q}(X)-\dim \ker \partial\bar{\partial}_{\phi(t)}\cap\left(\mathcal{H}_{BC}^{p,q}(X)+\Image(\partial\bar{\partial})^*\right)\cap A^{p,q}(X)\geq 0,
\]
then we have
\begin{equation}\label{eq-1461}
\dim H_{BC}^{p,q}(X)=\dim H_{BC\phi(t)}^{p,q}(X)+v^{p,q}_t+u^{p-1,q-1}_t.
\end{equation}
\end{theorem}
\begin{proof}This follows immediately from Lemma \ref{lem-ddbar*-closed-rep of deformed-BC-coho} and Proposition \ref{prop-kerddbar*-imageddbarphi}.
\end{proof}

\section{The deformed Fr\"ohlicher spectral sequences and the $\partial\bar{\partial}_{\phi}$-lemma}
\label{The deformed Frohlicher spectral sequences and ddbar-lemma}
Let $X$ be a complex manifold and $X_t$ a small deformation (of $X$) whose complex structure is represented by a Beltrami differential $\phi\in A^{0,1}(X, T_{X}^{1,0})$. Set the deformed de Rahm cohomology as
\[
H^{\bullet}_{d_\phi}(X):=\ker d_\phi/\Image d_\phi,
\]
then it is clear that $e^{i_{\phi}}:H^{\bullet}_{d_\phi}(X)\to H^{\bullet}_{dR}(X)$ is an isomorphism and the identity map induces the following commutative diagram:
\begin{equation}\label{dia-deofrmed 131}
\xymatrix{
 & H^{\bullet,\bullet}_{BC\phi}(X) \ar[ld] \ar[d] \ar[rd] & \\
 H^{\bullet,\bullet}_{\partial}(X) \ar[rd] & H^{\bullet,\bullet}_{d_\phi}(X) \ar[d] & H^{\bullet,\bullet}_{\bar{\partial}_\phi}(X) \ar[ld] \\
 & H^{\bullet,\bullet}_{A\phi}(X) & \\
}
\end{equation}
\begin{definition}\label{def-deformed Frohlicher spectral sequences and ddbar-lemma}
The spectral sequence associated to the double complex $(A^{\bullet,\bullet}(X), \partial, \bar{\partial}_{\phi})$ will be called the \emph{deformed Fr\"ohlicher spectral sequence} and we say $X$ satisfies the $\partial\bar{\partial}_{\phi}$-lemma if the homomorphism $H^{\bullet,\bullet}_{BC\phi}(X)\to H^{\bullet,\bullet}_{d_\phi}(X)$ in \eqref{dia-deofrmed 131} is injective, i.e.
\[
\ker \partial \cap \ker \bar{\partial}_{\phi} \cap \Image d_\phi = \Image\partial\bar{\partial}_{\phi}.
\]
\end{definition}
Set $d_\phi^c:=J^{-1}d_\phi J=\sqrt{-1}(\bar{\partial}_{\phi}-\partial)$, where $J$ is the almost complex structure on $X$. It is easy to see that $\ker\partial\cap\ker\bar{\partial}_{\phi}=\ker d_\phi\cap \ker d_\phi^c$ and $\Image\partial\bar{\partial}_{\phi}=\Image d_\phi d_\phi^c$. Hence, $X$ satisfies the $\partial\bar{\partial}_{\phi}$-lemma if and only if
\begin{equation}\label{eq-dphi-dphic-lemma-1}
\ker d_\phi \cap \ker d_\phi^c \cap \Image d_\phi = \Image d_\phi d_\phi^c~,
\end{equation}
or
\begin{equation}\label{eq-dphi-dphic-lemma-2}
\ker d_\phi \cap \ker d_\phi^c \cap \Image d_\phi^c = \Image d_\phi d_\phi^c~.
\end{equation}
There are two natural filtrations on $A^{\bullet,\bullet}(X)$:
\[
F^pA^{k}(X)=\bigoplus_{p\leq r\leq k}A^{r,k-r}(X),~\bar{F}^pA^{k}(X)=\bigoplus_{p\leq s\leq k}A^{k-s,s}(X),
\]
which induces two filtrations on the deformed de Rahm cohomology $H^{k}_{d_\phi}(X)$ for each $k\geq 0$:
\[
F^pH^{k}_{d_\phi}(X)=\{\alpha\in H^{k}_{d_\phi}(X)\mid \exists u\in F^pA^{k}(X)~s.t.~\alpha=[u] \},
\]
and
\[
\bar{F}^pH^{k}_{d_\phi}(X)=\{\alpha\in H^{k}_{d_\phi}(X)\mid \exists u\in \bar{F}^pA^{k}(X)~s.t.~\alpha=[u] \}.
\]
As usual, there are many ways to characterize the $\partial\bar{\partial}_{\phi}$-lemma:
\begin{proposition}\label{prop-char-ddbar-phi}The following statements are equivalent:
\begin{itemize}
  \item[1.]$X$ satisfies the $\partial\bar{\partial}_{\phi}$-lemma;
  \item[2.]The maps in \eqref{dia-deofrmed 131} induced by the identity map are all isomorphisms;
  \item[3.]The deformed Fr\"ohlicher spectral sequence degenerates at $E_1$ and there is a Hodge decomposition
\[
H^{k}_{d_\phi}(X;\mathbb{C})=\bigoplus_{p+q=k}F^pH^{k}_{d_\phi}(X)\cap \bar{F}^qH^{k}_{d_\phi}(X)~,~\forall k .
\]
\end{itemize}
\end{proposition}
\begin{proof}This follows directly from ~\cite[pp.\,268]{DGMS75}.
\end{proof}

\begin{theorem}\label{thm-AT-Frolicher-inequality-deformed}
Let $X$ be a compact complex manifold and $X_t$ a small deformation (of $X$) whose complex structure is represented by a Beltrami differential $\phi\in A^{0,1}(X, T_{X}^{1,0})$. Then for every $(p,q)\in \mathbb{N}\times \mathbb{N}$, we have
\begin{equation}\label{eq-AT-Frolicher-inequality-deformed}
\dim H_{BC\phi}^{p,q}(X)+ \dim H_{A\phi}^{p,q}(X) \geq \dim H_{\bar{\partial}_t}^{p,q}(X_t)+ \dim H_{\partial}^{p,q}(X).
\end{equation}
In particular, for every $k\in \mathbb{N}$, we have
\begin{equation}\label{eq-ddbar-phi-lemma-kehua}
\sum_{p+q=k}\dim H_{BC\phi}^{p,q}(X)+ \sum_{p+q=k}\dim H_{A\phi}^{p,q}(X) \geq 2\dim H_{dR}^{k}(X),
\end{equation}
and equality holds if and only if $X$ satisfies the $\partial\bar{\partial}_{\phi}$-lemma.
\end{theorem}
\begin{proof}This follows from similar arguments as in~\cite{AT13}. In fact, this theorem is a direct consequence of ~\cite[Thm.\,1\,and\,2]{AT15b} by noting that $\dim H_{\bar{\partial}_{\phi}}^{p,q}(X)=\dim H_{\bar{\partial}_t}^{p,q}(X_t)$~\cite[Thm.\,4.4]{Xia19dDol}.
\end{proof}
\begin{remark}\label{rk-to-AK-Frolicher-ineq}
\begin{itemize}
  \item[1.] From the work of Angella-Tardini~\cite[Thm.\,3.1]{AT17} we know that $X$ satisfies the $\partial\bar{\partial}_{\phi}$-lemma if and only if
\[
\sum_{p+q=k}\dim H_{BC\phi}^{p,q}(X)= \sum_{p+q=k}\dim H_{A\phi}^{p,q}(X)~;
\]
  \item[2.] From Proposition \ref{prop-char-ddbar-phi} we see that if $X$ satisfies the $\partial\bar{\partial}_{\phi}$-lemma, then for every $(p,q)\in \mathbb{N}\times \mathbb{N}$, we have
\[
\dim H_{BC\phi}^{p,q}(X)= \dim H_{A\phi}^{p,q}(X) = \dim H_{\bar{\partial}_t}^{p,q}(X_t)= \dim H_{\partial}^{p,q}(X).
\]
In particular, by Theorem \ref{thm-AT-Frolicher-inequality-deformed} we have $h_{BC\phi}^{k}=h_{A\phi}^{k}=h_{\bar{\partial}_t}^{k}(X_t)=h_{\bar{\partial}}^{k}=b_k$\footnote{We follow the notations as given in \cite{AT13}, e.g. $h_{BC\phi}^{k}:=\sum_{p+q=k}\dim H_{BC\phi}^{p,q}(X)$ and $b_k$ is the $k$-th Betti number.}, namely, the Fr\"ohlicher spectral sequence of $(A^{\bullet,\bullet}(X), \partial, \bar{\partial})$ degenerates at $E_1$.
\end{itemize}
\end{remark}

\begin{corollary}\label{coro-deformation-closedness-ddbar-lemma}
Let $\pi: (\mathcal{X}, X)\to (B,0)$ be a small deformation of the compact complex manifold $X$ such that for each $t\in B$ the complex structure on $X_t$ is represented by Beltrami differential $\phi(t)$. Then
the set
\[
T:=\{t\in B\mid X~\text{satisfies the}~\partial\bar{\partial}_{\phi(t)}\text{-lemma} \}
\]
is an analytic open subset (i.e. complement of analytic subset) of $B$. In particular, if $B\subset \mathbb{C}$ is a small open disc with $0\in B$ and $T$ is not empty, then $T=B$ or $T=B\setminus\{0\}$.
\end{corollary}
\begin{proof}First, by Theorem~\ref{thm-AT-Frolicher-inequality-deformed}, $X$ satisfies the $\partial\bar{\partial}_{\phi(t)}$-lemma if and only if
\begin{equation}\label{eq-ddbar-phi-lemma-kehua=}
h_{BC\phi(t)}^{k}+ h_{A\phi(t)}^{k} = 2b_k.
\end{equation}
We note that by Theorem~\ref{th-analytic subset} the set $\{t\in B\mid \eqref{eq-ddbar-phi-lemma-kehua=}~\text{holds}\}$ is an analytic open subset of $B$ since
\[
\{t\in B\mid h_{BC\phi(t)}^{k}+ h_{A\phi(t)}^{k} = 2b_k,~\forall~k\}= B\setminus \{t\in B\mid h_{BC\phi(t)}^{k}+ h_{A\phi(t)}^{k}\geq 2b_k+1,~\forall~k\}.
\]
\end{proof}
In particular, if $X$ satisfies the $\partial\bar{\partial}$-lemma, then by the above corollary $X$ also satisfies the $\partial\bar{\partial}_{\phi(t)}$-lemma for any small $t\in B$. Combining this with Remark~\ref{rk-to-AK-Frolicher-ineq}, we get that the Hodge numbers $\dim H_{\bar{\partial}_t}^{p,q}(X_t)$ and $\dim H_{BC\phi(t)}^{p,q}(X)$ are independent of $t$.

Recall that a smooth manifold $X$ is called \emph{formal} if its de Rahm complex $(A^{\bullet}(X),d)$ is formal as a differential graded algebra (DGA for short). The later means that there is a sequence of quasi-isomorphisms from $(A^{\bullet}(X),d)$ to its cohomology algebra $(H_{dR}^{\bullet}(X),0)$\footnote{Here, $(H^{\bullet}(X),0)$ is considered as a differential graded algebra with trivial differential.}, see~\cite{DGMS75,FHT01}.
\begin{theorem}\label{thm-formal}
Let $X$ be a compact complex manifold and $X_t$ a small deformation (of $X$) whose complex structure is represented by a Beltrami differential $\phi\in A^{0,1}(X, T_{X}^{1,0})$. If $X$ satisfies the $\partial\bar{\partial}_\phi$-lemma, then $X$ is formal.
\end{theorem}
\begin{proof}Consider the following homomorphisms of DGA
\begin{equation*}
\xymatrix@C=0.5cm{
(A^{\bullet}(X),d_\phi) && \ar[ll]_{i\quad}  (A^{\bullet}(X)\cap\ker d_\phi^c ,d_\phi)   \ar[rr]^{\quad p} && (H^{\bullet}_{d_\phi^c}(X),d_\phi=0),}
\end{equation*}
where $i$ is the inclusion and $p$ is the projection. We claim that the induced map $i^*$ is an isomorphism on cohomology. Indeed, $\forall x\in \ker d_\phi\cap \ker d_\phi^c$ if $x\in \Image d_\phi$ then by \eqref{eq-dphi-dphic-lemma-1}, $x\in\Image d_\phi d_\phi^c\Rightarrow i^*$ is injective; on the other hand, by \eqref{eq-dphi-dphic-lemma-2} $\forall x\in \ker d_\phi$ there exist $y\in A^{\bullet}(X)$ such that $x-d_\phi y\in \ker d_\phi\cap \ker d_\phi^c$, this shows that $i^*$ is surjective. Similarly, one shows that $p^*$ is an isomorphism on cohomology and $d_\phi=0$ on $H^{\bullet}_{d_\phi^c}(X)$. The conclusion then follows since $(A^{\bullet}(X),d_\phi)$ is isomorphic to $(A^{\bullet}(X),d)$ and $H^{\bullet}_{d_\phi^c}(X)\cong H^{\bullet}_{d_\phi}(X)\cong H^{\bullet}_{dR}(X)$.
\end{proof}

\section{The deformed Bott-Chern cohomology of the Iwasawa manifold and the holomorphically parallelizable Nakamura manifold}
\label{The deformed Bott-Chern cohomology of the Iwasawa manifold and the holomorphically parallelizable Nakamura manifold}
\begin{example}
Case III-(2). Let $G$ be the matrix Lie group defined by
\[
G := \left\{
\left(
\begin{array}{ccc}
 1 & z^1 & z^3 \\
 0 &  1  & z^2 \\
 0 &  0  &  1
\end{array}
\right) \in \mathrm{GL}(3;\mathbb{C}) \mid z^1,\,z^2,\,z^3 \in\mathbb{C} \right\}\cong \mathbb{C}^3
\;.
\]
Consider the discrete subgroup $\Gamma$ defined by
\[
\Gamma := \left\{
\left(
\begin{array}{ccc}
 1 & \omega^1 & \omega^3 \\
 0 &  1  & \omega^2 \\
 0 &  0  &  1
\end{array}
\right) \in G \mid \omega^1,\,\omega^2,\,\omega^3 \in\mathbb{Z}[\sqrt{-1}] \right\}
\;,
\]
The quotient $X=G/\Gamma$ is called the \emph{Iwasawa manifold}. A basis of $H^0(X,\Omega^1)$ is given by
\[
\varphi^1 = d z^1,~ \varphi^2 = d z^2,~ \varphi^3 = d z^3-z^1\,d z^2,
\]
and a dual basis $\theta^1, \theta^2, \theta^3\in H^0(X,T_X^{1,0})$ is given by
\[
\theta^1=\frac{\partial}{\partial z^1},~\theta^2=\frac{\partial}{\partial z^2} + z^1\frac{\partial}{\partial z^3},~\theta^3=\frac{\partial}{\partial z^3}.
\]
$X$ is equipped with the Hermitian metric $\sum_{i=1}^3\varphi^i\otimes\bar{\varphi}^i$. The Beltrami differential of the Kuranishi family of $X$ is
\[
\phi(t) = t_{i\lambda}\theta^i\bar{\varphi}^{\lambda} - D(t)\theta^3\bar{\varphi}^{3},~\text{with}~D(t)=t_{11}t_{22}-t_{21}t_{12},
\]
and the Kuranishi space of $X$ is
\[
\mathcal{B}=\{t=(t_{11}, t_{12}, t_{21}, t_{22}, t_{31}, t_{32})\in \mathbb{C}^6\mid |t_{i\lambda}|<\epsilon, i=1, 2, 3, \lambda=1,2 \},
\]
where $\epsilon>0$ is sufficiently small. Set
\[
\phi_1=\sum_{i=1}^3\sum_{\lambda=1}^2t_{i\lambda}\theta^i\bar{\varphi}^{\lambda},~ \phi_2 = D(t)\theta^3\bar{\varphi}^{3},
\]
and write the canonical deformation of $\sigma_0\in H_{BC}^{p,q}(X)$ by $\sigma(t)=\sum_k \sigma_k$ with each
\[
\sigma_k= - G_{BC}A\sum_{i+j=k} \partial i_{\phi_j} \sigma_i ,
\]
being the homogeneous term of degree $k>0$ in $t\in \mathcal{B}$. We will use the isomorphism $H_{BC\phi(t)}^{p,q}(X)\cong\dim V_t/\ker f_t$ proved in Proposition \ref{prop V-E-t} to compute $\dim H_{BC\phi(t)}^{p,q}(X)$. Since $\mathcal{B}$ is a polydisc, it is sufficient to check the coefficients of $d_{\phi(t)}\sigma(t)=0$, that is,
\begin{equation}\label{obstruction eq in example}
\partial\sigma_k= \bar{\partial}\sigma_k+ \sum_{j=1}^{k} \partial i_{\phi_j} \sigma_{k-j}=0 , ~~k>0~ .
\end{equation}
\end{example}

Let us now consider Bott-Chern deformations of forms in the harmonic space:
\[
\mathcal{H}_{BC}^{2,2}(X)=\mathbb{C}\{\varphi^{12\overline{13}}, \varphi^{12\overline{23}}, \varphi^{13\overline{12}}, \varphi^{13\overline{13}}, \varphi^{13\overline{23}}, \varphi^{23\overline{12}}, \varphi^{23\overline{13}}, \varphi^{23\overline{23}}\}.
\]
Set $\sigma_0 = \sum a_{ijkl}\varphi^{ij\overline{kl}}\in \mathcal{H}_{BC}^{2,2}(X)$, then
\[
\partial i_{\phi_1}\sigma_0= (-t_{12}a_{13\overline{13}} + t_{11}a_{13\overline{23}} - t_{22}a_{23\overline{13}} + t_{21}a_{23\overline{23}}) \varphi^{12\overline{123}}
\]
is $\bar{\partial}$-exact if and only if
\begin{equation}\label{a_{23}}
t_{12}a_{13\overline{13}} - t_{11}a_{13\overline{23}} + t_{22}a_{23\overline{13}} - t_{21}a_{23\overline{23}}=0 ,
\end{equation}
and in this case
\[
\sigma_{1}=-G_{BC}A \partial i_{\phi_1}\sigma_0 =0~ .
\]
But
\[
\partial i_{\phi_2}\sigma_0 = 0 \Longrightarrow \sigma_{2}=-G_{BC}A \partial  ( i_{\phi_2}\sigma_0 + i_{\phi_1}\sigma_1 ) =0~ ,
\]
and $\phi_k=0,~ k>2$ we thus have $\sigma_{k}=0,~ k>2$.

Therefore, for $V=\mathcal{H}_{BC}^{2,2}(X)$ we have (see Definition \ref{def-V_t})
\begin{align*}
V_{t}=
&\{ \sum a_{ijkl}\varphi^{ij\overline{kl}}\in \mathcal{H}_{BC}^{2,2}(X) \mid  (a_{1213}, a_{1223}, a_{1312}, a_{1313}, a_{1323}, a_{2312}, a_{2313}, a_{2323})\in \mathbb{C}^8\\&~\text{s.t.}~\sigma(t)\in\ker d_{\phi(t)},
\text{where}~\sigma (t)=\sum_{k} \sigma_k~\text{with}~ \sigma_0=\sum a_{ijkl}\varphi^{ij\overline{kl}}\\
&~\text{and}~ \sigma_{k}=-G_{BC}A\sum_{i+j=k} \partial i_{\phi_j} \sigma_i,~\forall k\neq 0 \}\\
=&\{\sum a_{ijkl}\varphi^{ij\overline{kl}} \mid (a_{1213}, a_{1223}, a_{1312}, a_{1313}, a_{1323}, a_{2312}, a_{2313}, a_{2323})\in \mathbb{C}^8~\text{satisfy}~ \eqref{a_{23}}\}.
\end{align*}
On the other hand, $\Image\partial\bar{\partial}_{\phi(t)}=\mathbb{C}\{\partial\bar{\partial}_{\phi(t)}\varphi^{3\bar{3}}\}=\mathbb{C}\{\varphi^{12\overline{12}}\}$ (Since $X$ is parallelizable, we only need to consider left invariant forms. See the discussions in the last paragraph of this section) and
\[
(\partial\bar{\partial})^*\varphi^{12\overline{12}}=(*\partial*)(*\bar{\partial}*)\varphi^{12\overline{12}}
=-(*\partial\bar{\partial}*)\varphi^{12\overline{12}}=-*\partial\bar{\partial}\varphi^{3\bar{3}}=\varphi^{3\bar{3}}\neq 0,
\]
implies
\[
\ker f_t\cong \ker(\partial\bar{\partial})^*\cap\Image\partial\bar{\partial}_{\phi(t)} \cap A^{2,2}(X)=0.
\]
By Proposition \ref{prop V-E-t} we have
\begin{equation}
\dim H_{BC\phi(t)}^{2,2}(X)=\dim V_t-\dim\ker f_t=
\left\{
\begin{array}{rcl}
8,~ &(t_{11}, t_{12}, t_{21}, t_{22})=0 \\[5pt]
7,~ &(t_{11}, t_{12}, t_{21}, t_{22})\neq0 ~.
\end{array}
\right.
\end{equation}

The other deformed Bott-Chern cohomology can be computed in the same way. Write $(i), (ii), (iii)$ for the three cases when $(t_{11}, t_{12}, t_{21}, t_{22})=0$, $(t_{11}, t_{12}, t_{21}, t_{22})\neq0$ and $D(t)= 0$, $D(t)\neq 0$, respectively. Then we have the following (where $h^{p,q}:=\dim H_{BC\phi(t)}^{p,q}(X)$ and $t\in (i), (ii), (iii)$, respectively)
\renewcommand\arraystretch{1.5}
\begin{table}[!htbp]
\centering
\begin{center}
\begin{tabular}{|cc|ccc|cccc|ccc|cc|}
\hline
 $h^{1,0}$ & $h^{0,1}$ & $h^{2,0}$ & $h^{1,1}$ & $h^{0,2}$ & $h^{3,0}$ & $h^{2,1}$ & $h^{1,2}$ & $h^{0,3}$ & $h^{3,1}$ & $h^{2,2}$ & $h^{1,3}$ & $h^{3,2}$ & $h^{2,3}$\\
\hline
 $2$ & $2$ & $3$ & $4$ & $3$ & $1$ & $6$ & $6$ & $1$ & $2$ & $8$ & $2$ & $3$ & $3$ \\
\hline
 $2$ & $2$ & $2$ & $4$ & $3$ & $1$ & $6$ & $6$ & $1$ & $2$ & $7$ & $2$ & $3$ & $3$\\
\hline
 $2$ & $2$ & $1$ & $4$ & $3$ & $1$ & $6$ & $6$ & $1$ & $2$ & $7$ & $2$ & $3$ & $3$\\
\hline
\end{tabular}
\end{center}
\end{table}

Comparing this with the computations made by Angella~\cite{Ang13} we see that $\dim H_{BC\phi(t)}^{p,q}(X)=\dim H_{BC}^{p,q}(X_t)$ is not true in general for $p=q$.

\begin{example}
Case III-(3b). Let $X=\mathbb{C}^3/\Gamma$ be the solvable manifold constructed by Nakamura in Example III-(3b) of~\cite{Nak75}. We have
\begin{align*}
H^0(X,\Omega_X^{1})~ =& ~\mathbb{C}\{\varphi^1 = d z^1,~ \varphi^2 = e^{z_1} d z^2,~ \varphi^3 = e^{-z_1} d z^3\} ~,\\
H^0(X,T_X^{1,0})~  =& ~\mathbb{C}\{\theta^1=\frac{\partial}{\partial z^1},~\theta^2= e^{-z_1} \frac{\partial}{\partial z^2},~\theta^3= e^{z_1} \frac{\partial}{\partial z^3}\} ~,\\
\mathcal{H}^{0,1}(X)~ =& ~\mathbb{C}\{ \psi^{\bar{1}} = d z^{\bar{1}},~ \psi^{\bar{2}} = e^{z_1} d z^{\bar{2}},~ \psi^{\bar{3}} = e^{-z_1} d z^{\bar{3}}\} ~,\\
\mathcal{H}^{0,1}(X, T_X^{1,0})~ =& ~\mathbb{C}\{\theta^i\psi^{\bar{\lambda}} ,~ i=1, 2, 3, \lambda=1, 2, 3\} ~,
\end{align*}
where $X$ is equipped with the Hermitian metric $\sum_{i=1}^3\varphi^i\otimes\bar{\varphi}^i$. The Beltrami differential of the Kuranishi family of $X$ is
\[
\phi(t) = \phi_1 = t_{i\lambda}\theta^i\psi^{\bar{\lambda}}
\]
and the Kuranishi space of $X$ is
\[
\mathcal{B}=\{t=(t_{11}, t_{12}, t_{13}, t_{21}, t_{22}, t_{23}, t_{31}, t_{32}, t_{33})\in \mathbb{C}^9\mid |t_{i\lambda}|<\epsilon, i=1, 2, 3, \lambda=1, 2, 3 \},
\]
where $\epsilon>0$ is sufficiently small. We will restrict to the one parameter family defined by $t_{12}=t_{13}=t_{21}=t_{22}=t_{23}=t_{31}=t_{32}=t_{33}=0$ and in this case the Beltrami differential is $\phi=\phi(t)=t\frac{\partial}{\partial z^1}dz^{\bar{1}}$ where $t=t_{11}$.

Let us consider the Bott-Chern deformations of forms in
\begin{align*}
\mathcal{H}_{BC}^{2,1}(X)=\mathbb{C}\{&e^{z^1}dz^{12\overline{1}}, e^{2z^1}dz^{12\overline{2}}, dz^{12\overline{3}}, e^{-z^1}dz^{13\overline{1}}, dz^{13\overline{2}}, e^{-2z^1}dz^{13\overline{3}}, dz^{23\overline{1}},\\ &e^{z^{\bar{1}}}dz^{13\overline{1}}, e^{z^{\bar{1}}}dz^{12\overline{1}}\}.
\end{align*}
Set
\begin{align*}
\sigma_0 =& a_{121}e^{z^1}dz^{12\overline{1}}+a_{122}e^{2z^1}dz^{12\overline{2}}+a_{123}dz^{12\overline{3}}+a_{131}e^{-z^1}dz^{13\overline{1}}+a_{132}dz^{13\overline{2}}\\&+a_{133}e^{-2z^1}dz^{13\overline{3}} +a_{231}dz^{23\overline{1}}+b_{131}e^{z^{\bar{1}}}dz^{13\overline{1}}+b_{121}e^{z^{\bar{1}}}dz^{12\overline{1}},
\end{align*}
then
\[
\partial i_{\phi_1}\sigma_0= -2a_{122}te^{2z^1}dz^{12\overline{12}} + 2a_{133}te^{-2z^1}dz^{13\overline{13}}
\]
is $\bar{\partial}$-exact if and only if $t=0$.
Therefore, for $V=\mathcal{H}_{BC}^{2,1}(X)$ and $t\neq 0$ we have
\begin{align*}
V_{t}=
&\{ \sigma_0\in \mathcal{H}_{BC}^{2,1}(X) \mid  (a_{121}, a_{122}, a_{123}, a_{131}, a_{132}, a_{133}, a_{231}, b_{131}, b_{121} )\in \mathbb{C}^9~\text{s.t.}~\\&\sigma(t)\in\ker d_{\phi(t)},
\text{where}~\sigma (t)=\sum_{k} \sigma_k~\text{with}~\sigma_{k}=-G_{BC}A\sum_{i+j=k} \partial i_{\phi_j} \sigma_i,~\forall k\neq 0 \}\\
=&\mathbb{C}\{e^{z^1}dz^{12\overline{1}}, dz^{12\overline{3}}, e^{-z^1}dz^{13\overline{1}}, dz^{13\overline{2}}, dz^{23\overline{1}}, e^{z^{\bar{1}}}dz^{13\overline{1}}, e^{z^{\bar{1}}}dz^{12\overline{1}}\}.
\end{align*}
On the other hand,
\begin{align*}
\ker(\partial\bar{\partial})^*\cap\Image\partial\bar{\partial}_{\phi(t)} \cap A^{2,1}(X)=&\mathbb{C}\{\partial\bar{\partial}_{\phi(t)}e^{z^1}dz^{2}, \partial\bar{\partial}_{\phi(t)}e^{-z^1}dz^{3}\}\\=&\mathbb{C}\{te^{z^1}dz^{12\bar{1}}, te^{-z^1}dz^{13\bar{1}}\}
\end{align*}
and
\[
\dim\ker f_t= \dim \ker(\partial\bar{\partial})^*\cap\Image\partial\bar{\partial}_{\phi(t)} \cap A^{2,2}(X)=2 .
\]
By Proposition \ref{prop V-E-t} we have
\begin{equation}
\dim H_{BC\phi(t)}^{2,1}(X)=\dim V_t-\dim\ker f_t=
\left\{
\begin{array}{rcl}
9,~ &t=0 \\[5pt]
5,~ &t\neq0 ~.
\end{array}
\right.
\end{equation}
We summarise the computations of the deformed Bott-Chern cohomology in this case as follows (where $h^{p,q}:=\dim H_{BC\phi(t)}^{p,q}(X)$ and $t=0,~\neq 0$, respectively):
\renewcommand\arraystretch{1.5}
\begin{table}[!htbp]
\centering
\begin{center}
\begin{tabular}{|cc|ccc|cccc|ccc|cc|}
\hline
 $h^{1,0}$ & $h^{0,1}$ & $h^{2,0}$ & $h^{1,1}$ & $h^{0,2}$ & $h^{3,0}$ & $h^{2,1}$ & $h^{1,2}$ & $h^{0,3}$ & $h^{3,1}$ & $h^{2,2}$ & $h^{1,3}$ & $h^{3,2}$ & $h^{2,3}$\\
\hline
 $1$ & $1$ & $3$ & $7$ & $3$ & $1$ & $9$ & $9$ & $1$ & $3$ & $11$ & $3$ & $5$ & $5$ \\
\hline
 $1$ & $1$ & $1$ & $5$ & $3$ & $1$ & $5$ & $7$ & $1$ & $1$ & $7$  & $3$ & $3$ & $3$\\
\hline
\end{tabular}
\end{center}
\end{table}

From this table and~\cite{AK17b}, we notice that $X_t$ satisfy the $\partial\bar{\partial}$-lemma but $X$ does not satisfy the $\partial\bar{\partial}_{\phi(t)}$-lemma for any $t\neq 0$ .
\end{example}

We need to point out that in the above computations (especially those concerning $\ker(\partial\bar{\partial})^*\cap\Image\partial\bar{\partial}_{\phi(t)}$), only invariant forms is considered. This is valid because the Bott-Chern cohomology of complex parallelizable manifold may be computed by left invariant forms~\cite{Ang13} and given a family of deformations $\{X_t\}_{t\in B}$ of such manifolds the set of $t$ for which the deformed Bott-Chern cohomology may be computed by left invariant forms is an open subset of $B$ (this will be proved in \cite{Xia20g}).


\vskip 1\baselineskip \textbf{Acknowledgements.} I would like to thank Prof. Kefeng Liu for his constant encouragement and many useful discussions. Many thanks to Shengmao Zhu, Sheng Rao and Daniele Angella for useful communications. I would also like to thank Prof. Bing-Long Chen for his constant support.

\bibliographystyle{alpha}

\end{document}